\documentclass{emsprocart}
\usepackage{graphics, makeidx}
 \usepackage[colorlinks=true, pdfstartview=FitV, linkcolor=blue, citecolor=blue, urlcolor=blue, breaklinks=true]{hyperref}

\newtheorem{theorem}{Theorem}[section]
\newtheorem{corollary}[theorem]{Corollary}

\newtheorem{proposition}[theorem]{Proposition}

\newtheorem{exam}[theorem]{Example}
\newtheorem{example}[theorem]{Example}
\newtheorem*{ack}{Acknowledgment}

\theoremstyle{definition}
\newtheorem{definition}[theorem]{Definition}
\newtheorem{remark}[theorem]{Remark}

\newcommand{\Hom}{\operatorname{Hom}}

\title[On toric degenerations of flag varieties]{On toric degenerations of flag varieties}

\author[Xin Fang, Ghislain Fourier and Peter Littelmann]{Xin Fang, Ghislain Fourier and Peter Littelmann}
\contact[xfang@math.uni-koeln.de]{Xin Fang, University of Cologne}
\contact[Ghislain.Fourier@glasgow.ac.uk]{Ghislain Fourier, University of Glasgow}
\contact[peter.littelmann@math.uni-koeln.de]{Peter Littelmann, University of Cologne}
\makeindex

\begin{document}

\begin{abstract}
Following the historical track in pursuing $T$-equivariant flat toric degenerations of flag varieties and spherical varieties, we explain how powerful tools in algebraic geometry and representation theory, such as canonical bases, Newton-Okounkov bodies, PBW-filtrations and cluster algebras
come to push the subject forward.
\end{abstract}



\maketitle
\setcounter{tocdepth}{1}
\tableofcontents
\section{Introduction}
One of the beautiful and astonishing properties of the theory of toric varieties is the powerful dictionary
translating algebraic geometry properties into combinatorial properties in terms of lattices, cones and polytopes, and vice versa.
It is hence tempting to extend this powerful machinery to a larger class of varieties by studying flat toric degenerations
of a given variety. 

Several recent developments have enhanced the interest in flat toric degenerations. One feature
has been the evolution of the theory of Newton-Okounkov bodies \cite{K1, KK} and  
its applications in algebraic and in symplectic geometry, see, for example \cite{FLP, HP, HK,  K1, K2, KK, KüL,LaMu, NNU2, WN}. The theory of Newton-Okounkov
bodies attaches a monoid $\Gamma=\Gamma(X,\mathcal L)$ and a convex body $\Delta_{X,\mathcal L}$ 
to a polarized (smooth) complex projective variety $(X,\mathcal L)$.
If $\Gamma$ is finitely generated, then Anderson \cite{A}, Harada and Kaveh \cite{HK} show that there exists an integrable system $\mathcal F$ on $X$ 
such that the associated momentum map $\mu_{\mathcal F} : X \rightarrow \mathbb R^N$ has $\Delta_{X,\mathcal L}$ as image, 
and there exists a flat degeneration of $X$ into a projective toric variety $X_{\mathbb T}$ for a complex torus 
$\mathbb T\simeq (\mathbb C^*)^N$ having $\Delta_{X,\mathcal L}$ as associated momentum polytope. In addition,
the integrable system generates a torus action on the inverse image under $\mu_{\mathcal F}$ 
of the interior of the moment polytope, inducing a $\mathbb T$-equivariant symplectomorphism between this open subset  
of $X$ and an open subset of $X_{\mathbb T}$.

A second feature is the program of Gross, Hacking, Keel and Kontsevich \cite{GHKK}
concerning (partially compactified) cluster varieties (like flag varieties, double Bruhat cells etc.).
For a given cluster variety ${\cal X}$, let ${\cal X}^\vee$ be its mirror dual and denote
by  ${\cal X}^\vee({\mathbb R}^T)$ its tropicalization. For every seed one gets an identification 
of ${\cal X}^\vee({\mathbb R}^T)$ with a lattice in some $\mathbb R^{\dim \cal X}$.
This allows them to generalize basic polyhedral constructions from toric geometry
in a straight-forward way. In this setting, a polytope $\Xi\subset {\cal X}^\vee({\mathbb R}^T)$ (in the generalized sense) 
corresponds to a compactification $\widehat{\cal X}$ of ${\cal X}$. For any seed $S$ one gets 
an identification of the polytope $\Xi\subset {\cal X}^\vee({\mathbb R}^T)$ with a polytope 
$\Xi_S\subset {\mathbb R}^{\dim \cal X}$ and a flat degeneration 
of the compactification $\widehat{\cal X}$ to an ordinary polarized projective toric variety
with $\Xi_S$ as associated polytope.
In the program proposed by Gross, Hacking, Keel and Kontsevich, there are still many open questions
left. For example, in the case of flag varieties it is not clear yet whether the assumptions necessary 
for the program to work are fulfilled. Only recently Magee \cite{M} announced 
a proof of several conjectures on the cluster structure of $SL_n/U$. Gross, Hacking, Keel and Kontsevich 
conjecture \cite{GHKK} that their approach specializes to a uniform construction of many degenerations of representation 
theoretic objects to toric varieties. The results of Rietsch and Williams \cite{RW} suggest that for flag varieties
the two features: Newton-Okounkov theory and the cluster approach, are connected.

To determine explicitly the Newton-Okounkov body is in general a very difficult task. But to use the body to get,
for example, estimates for the Gromov width, the Seshadri constant, for symplectic packings  \cite{FLP, HP,  I, K2}, to determine potential functions 
\cite{NNU1,NNU2} etc.,
often a precise knowledge of the body is necessary. From this point of view the flag variety case can be seen
as a perfect toy model for the theory. Since it is endowed with a wealth of additional structure
(combinatorics of root systems and Weyl groups, representation theory, enveloping algebra, etc.), these extra features make it often 
possible to give an alternative construction of the monoid $\Gamma=\Gamma(X,\mathcal L)$ and to give 
a precise description of the convex body $\Delta_{X,\mathcal L}$. 
For this reason we will stick in this overview to the case of flag varieties and make sometimes
remarks about spherical varieties. For other connections of Newton-Okounkov bodies to representation theory see
also the work of H. Sepp\"anen, for example in \cite{Se}. 

We will not be able to discuss all degenerations of flag varieties which can be found in the literature
but just present a selection. For example, 
we will not be able to discuss the Gr\"obner degenerations of Gra\ss mann 
varieties arising from tropical geometry. More details on this can be found in \cite{SS}; 
a relation to the approach by Rietsch and Williams 
in \cite{RW} mentioned above can be found in \cite{BFFHL}.

We explain now in more detail the content of the various sections.
\par
In section ~\ref{Sec:Grassmannian} we will give a short introduction into the theory of Hodge algebras, also called algebras 
with straightening laws.
One of the first flat degenerations of a Gra\ss mann variety into a union of toric varieties has been constructed by Hodge \cite{Hodge}.
The ideas of Hodge have been largely generalized by De Concini, Eisenbud and Procesi in the framework of Hodge algebras 
\cite{DEP}.
\par
In section~\ref{Section:Chiviri} we consider an extended version of Hodge algebras.
For generalized flag varieties,  results similar to those of De Concini, Eisenbud and Procesi 
have been obtained by Chiriv\`i \cite{Ch1} via upgrading the Hodge algebra to the Lakshmibai-Seshadri (LS) algebra.
\par
In section~\ref{GnachU} we discuss the flat toric degenerations of flag varieties associated to string parametrizations.
A flat toric degeneration of the flag variety with an irreducible special fiber has been 
given by Gonciulea and Lakshmibai \cite{GL} in the case $SL_n/B$, 
where $B$ is a Borel subgroup, using standard monomial theory. This has been interpreted geometrically by Kogan and Miller \cite{KM} 
using geometric invariant theory. A uniform construction for arbitrary reductive algebraic groups $G$ has been given
by Caldero \cite{Ca1}. Let $G$ be a reductive algebraic group $G$ and fix a maximal unipotent subgroup $U$.
For every choice of a reduced decomposition $\underline{w_0}$ of the longest word $w_0$ in the Weyl group of $G$,
Caldero constructs a flat toric degeneration of the affine variety $G/\hskip -3.5pt/U$. 
The construction relies on the associated string parametrization and the multiplicative properties of the dual canonical bases. 
These results have been generalized to spherical varieties by Alexeev and Brion \cite{AB}, see also the articles by Kaveh
and Kiritchenko \cite{K1, Ki} for another 
approach via the framework of Newton-Okounkov bodies \cite{KK}.
\par
In section~\ref{Section:birational} we describe a new approach to determine  
Newton-Okounkov bodies for flag varieties \cite{FaFL}. Roughly speaking, the idea
is to replace $G/B$ by a birationally equivalent product of root subgroups.
The latter leads then naturally to a coordinate system for the function field of $G/B$
and an associated lowest term valuation. So each such product decomposition
gives rise to a valuation monoid $\Gamma_1$. The (up to birational equivalence) 
decomposition of $G/B$ as a product of root subgroups induces also a filtration
of the enveloping algebra $U(\mathfrak n^-)$ and of the finite dimensional irreducible representations
of $G$. This can be used to attach a second monoid $\Gamma_2$ to the decomposition, 
the monoid of essential monomials. It turns out that $\Gamma_1=\Gamma_2$, thus giving a representation theoretic 
interpretation of the valuation monoid. Once one knows that the monoid $\Gamma=\Gamma_1=\Gamma_2$
is finitely generated, the methods used by Alexeev and Brion \cite{AB} apply also to this much more general situation.
In particular, we can use the theory to describe degenerations of $G$-varieties and not only of flag varieties.

The background of the filtration construction is 
a conjecture of Vinberg, which leads to a new class of toric degenerations with irreducible special fibers,
see \cite{FFL1} for $G=SL_n$ and $Sp_{2n}$, see also \cite{G,G2} for $G$ of type $D_4$ and $G_2$.
This construction uses monomial bases obtained through a refinement of the PBW-filtration of the corresponding universal 
enveloping algebra \cite{FFL2,FFL3}.

In section~\ref{examplesectionASM} we discuss several examples. In particular, we show that the  combination of 
the approach via Newton-Okounkov bodies and the approach via filtrations provides a general framework to study toric degenerations of 
generalized flag varieties and spherical varieties: we recover the degenerations of flag varieties associated to string parametrizations
(section \ref{Sec:String}) as well as the toric degenerations induced by a PBW-filtration (section~\ref{GoodOrdering}),
and provide new examples arising from Lusztig parametrizations (section~\ref{Lusztig}) of the canonical basis. 

In section~\ref{Gromov} we present an application of the method and determine the Gromov width 
of a coadjoint orbit \cite{FLP}. 

In section~\ref{sec:smallrank} we discuss some small rank examples.
\par
In section~\ref{cluster} we give a very rough sketch on relations between the theory
of cluster varieties and toric degenerations \cite{GHKK,M,RW}.
To go into all the technical details would blow up the framework of this overview.
\par
We conclude the survey in section~\ref{openquestion} with some open questions.

\begin{ack}\rm The authors would like to thank St\'ephanie Cupit-Foutou, Konstanze Rietsch, 
Bea Schumann and Lauren Williams for helpful discussions and suggestions.

The first author is funded by the Alexander von Humboldt
Foundation. The research of the second and third authors has been partially supported
by the DFG-Priority Program SPP 1388 - Representation Theory.
\end{ack}
\section{Gra\ss mann varieties and algebras with straightening laws}\label{Sec:Grassmannian}
Flat degenerations of Gra\ss mann varieties have been already used by
Hodge \cite{Hodge}. The ideas of Hodge have been formalized and largely generalized later by
De Concini, Eisenbud and Procesi. For more details we recommend \cite{DEP}, where they introduce
the notion of a \textit{Hodge algebra}, also called an \textit{algebra with straightening laws}.
To make the presentation easier, we state the results only over the complex numbers.

Consider the Gra\ss mann variety $Gr_{d}(\mathbb C^n)\subseteq \mathbb P(\Lambda^d\mathbb C^n)$\index{$Gr_{d}(\mathbb C^n)$, Gra\ss mann variety} of
$d$-dimensional subspaces of $\mathbb C^n$. We view $Gr_{d}(\mathbb C^n)$ as a projective
subvariety of $\mathbb P(\Lambda^d\mathbb C^n)$ via the Pl\"ucker embedding, that is the embedding defined by sending a basis $(v_1,\ldots,v_d)$ of a $d$-dimensional
subspace $U$ to the class $[v_1\wedge\ldots\wedge v_d]\in \mathbb P(\Lambda^d\mathbb C^n)$.
By viewing matrices as sequences of column vectors, we can view this embedding as coming from the natural map
\begin{equation}\label{exteriorproduct}
\begin{matrix}
\pi_{d}:M_{n,d}(\mathbb C)&\to &\Lambda^{d}V\\
A=(v_1,\dots,v_d)&\mapsto&v_1\wedge \dots \wedge v_d\\
\end{matrix}.
\end{equation}
Note that the map is invariant with respect to the right multiplication by $SL_d(\mathbb C)$ on $M_{n,d}(\mathbb C)$.
It provides by classical invariant theory an identification of the categorical
quotient $M_{n,d}(\mathbb C)/\hskip-3.5pt/ SL_d(\mathbb C)$ with the affine cone $\widehat Gr_{d}(\mathbb C^n)$
over the Gra\ss mann variety. In particular, after replacing $M_{n,d}(\mathbb C)$ by the open subset
${M}_{n,d}(\mathbb C)-\pi_d^{-1}(0)$, we get the desired identifications: subspaces of dimension $d$ of $\mathbb C^n$
correspond to (by fixing a basis) $GL_n(\mathbb C)$-orbits in $M_{n,d}(\mathbb C)-\pi_d^{-1}(0)$, and these orbits
are (by the map above) in bijection with the points in $Gr_{d}(\mathbb C^n)$.

Let $I_{d,n}:=\{\underline{\mathbf i}=(i_1,\dots,i_d)\ |\ 1\le i_1<\dots<i_d\le n\}$ be the set of all strictly
increasing sequences of length $d$ between $1$ and $n$. The $d$-fold wedge product is alternating, the ordered
products of elements in the canonical basis of $\mathbb C^n$ form a basis
$\{e_{\underline{\mathbf i}}=e_{i_{1}}\wedge \cdots\wedge e_{i_{d}}$,  $\underline{\mathbf i}\in I_{d,n}\}$
of $\Lambda^d\,\mathbb C^n$, called the canonical basis of $\Lambda^d\,\mathbb C^n$.
Denote by $\{\mathbf p_{\underline{\mathbf i}}\mid \underline{\mathbf i}\in I_{d,n}\}$ its dual basis in $(\Lambda^{d} \mathbb C^n)^*$,
the $\mathbf p_{\underline{\mathbf i}}$\index{$\mathbf p_{\underline{\mathbf i}}$, Pl\"ucker coordinate} are called {\it Pl\"ucker coordinates}.

Denote by $\mathbb C[Gr_{d}(\mathbb C^n)]$ the homogeneous coordinate ring of the embedded
variety, it is the quotient of $\mathbb C[p_{\underline{\mathbf i}}\mid \underline{\mathbf i}\in I_{d,n}]$
by the vanishing ideal $I(Gr_{d}(\mathbb C^n))\subseteq \mathbb C[\Lambda^d \mathbb C^n]$
of the embedded Gra\ss mann variety $Gr_{d}(\mathbb C^n)\subseteq {\mathbb P}(\Lambda^d \mathbb C^n)$.
This ideal is generated by the {\it Pl\"ucker relations} (see \cite{BL,S}).

We define a partial order ``$\ge$'' on $I_{d,n}$ as follows:
$\underline{\mathbf i}\ge \underline{\mathbf j} \Leftrightarrow i_t\ge j_t$ for all $t=1,\ldots,d$.
A monomial $\mathbf p_{\underline{\mathbf i}^1}\cdots \mathbf p_{\underline{\mathbf i}^r}$ of Pl\"ucker coordinates
is called standard in $\mathbb C[\Lambda^d\,\mathbb C^n]$ iff ${\mathbf i}^1\ge \ldots \ge \underline{\mathbf i}^r$.
So what about monomials which are not standard?

The precise description of an algebra with straightening laws is somewhat technical,
we try to explain what happens in this special case. The second fundamental theorem in invariant theory
describes the relations among the Pl\"ucker coordinates
$\mathbf p_{\underline{\mathbf i}}$ considered as functions on ${\widehat Gr_{d}(\mathbb C^n)}$. If $\underline{\mathbf i}$ and $\underline{\mathbf j}$
are not comparable, then there exists a quadratic polynomial in $I(Gr_{d}(\mathbb C^n))$ such that
\begin{equation}\label{pluecker}
\mathbf p_{\underline{\mathbf i}}\mathbf p_{\underline{\mathbf j}}- \left(\mathbf p_{\underline{\mathbf i}\cup \underline{\mathbf j}} \mathbf p_{\underline{\mathbf i}\cap \underline{\mathbf j}} +
\hbox{other quadratic terms} \right) \in I(Gr_{d}(\mathbb C^n)).
\end{equation}
Here $\underline{\mathbf i}\cup \underline{\mathbf j}=(\max(i_1,j_1),\ldots,\max(i_d,j_d))$
and $\underline{\mathbf i}\cap \underline{\mathbf j}=(\min(i_1,j_1),\ldots,\min(i_d,j_d))$.

Refine the partial order "$\ge$" to a total order "$\succ$" and denote with the same symbol
the induced lexicographic order on the polynomial ring $\mathbb C[\Lambda^d\,\mathbb C^n]$,
which is a monomial order. One can show: all the terms in brackets of \eqref{pluecker} are strictly
larger than $\mathbf p_{\underline{\mathbf i}}\mathbf p_{\underline{\mathbf j}}$.

So if one has still some non-standard monomials in the bracket of \eqref{pluecker},
one may repeat the procedure with these terms again. Since the new terms occurring are always of degree two
and strictly larger in the monomial order, after a finite number of repetitions one sees: a monomial, which is not standard,
can be expressed modulo $I(Gr_{d}(\mathbb C^n))$ as a sum of standard monomials. In addition, all these standard monomials are larger in the
monomial order than the monomial we started with. Such an algorithm expressing non standard monomials as a linear combination of
standard monomials is called a {\it straightening law}.

An algebra with such properties: a basis consisting of a special class of monomials, the {\it standard monomials},
together with relations expressing non standard monomials as a linear combination of
(larger) standard monomials, this is roughly what is called an {\it algebra with straightening laws}.

In the case of the Gra\ss mann variety we get more precisely (for more details and other fields see \cite{BL},  \cite{DEP} and \cite{S}):
\begin{theorem}\label{grassmanncase}
\begin{enumerate}
\item[{\it i)}] The standard monomials form a basis of $\mathbb C[Gr_{d}(\mathbb C^n)]$.
\item[{\it ii)}] $I(Gr_{d}(\mathbb C^n))$ is generated by the straightening laws.
\item[{\it iii)}] There exists a flat degeneration of $Gr_{d}(\mathbb C^n)$ into a union of toric varieties,
taking the straightening relations to their initial terms. I.e., the vanishing ideal
of the initial scheme is generated by the monomials $\mathbf p_{\underline{\mathbf i}}\mathbf p_{\underline{\mathbf j}}$
for all pairs $(\underline{\mathbf i}, \underline{\mathbf j})$ such that $\underline{\mathbf i}$ and $\underline{\mathbf j}$
are not comparable with respect to the partial order "$\ge$". The initial scheme is a union of projective spaces,
one for each maximal chain in the partially ordered set $I_{d,n}$.
\end{enumerate}
\end{theorem}
\par\noindent
\begin{exam}\rm
The Gra\ss mann variety $Gr_{2}(\mathbb C^4)\subset \mathbb P(\Lambda^2\mathbb C^4)$ is defined as the zero set
of the homogeneous equation
$$
\mathbf p_{\mathbf{[14]}} \mathbf p_\mathbf{[23]} -
\mathbf p_{\mathbf{[24]}}\mathbf p_{\mathbf{[13]}}+
\mathbf p_{\mathbf{[34]}}\mathbf p_{\mathbf{[12]}}=0.
$$
The partial order on the set $I_{2,4}$ is the following:
$$
\begin{array}{ccccccc}
\mathbf{[12]}&<& \mathbf{[13]}&<& \mathbf{[14]}\\
                     &  &     \rotatebox{90}{$>$}   &&   \rotatebox{90}{$>$} \\
                     &   & \mathbf{[23]}&<&\mathbf{[24]}&<&\mathbf{[34]}.
\end{array}
$$
The Pl\"ucker relation above expresses the only quadratic non-standard monomial $\mathbf p_{\mathbf{[14]}} \mathbf p_\mathbf{[23]}$
as a linear combination of two standard monomials, modulo the vanishing ideal.
The degenerate version of $Gr_{2}(\mathbb C^4)$ is the reducible variety defined as the zero set of the homogeneous equation
$\mathbf p_{\mathbf{[14]}} \mathbf p_\mathbf{[23]}=0$.
\end{exam}

\section{Chiriv\`i's extension to the generalized flag varieties}\label{Section:Chiviri}
Let $G$\index{$G$, complex algebraic group} be a complex semisimple algebraic group, we fix a Borel subgroup $B$\index{$B$, Borel subgroup} and a maximal torus $T$\index{$T$, torus}. Let $\Lambda$\index{$\Lambda$, weight lattice} be the weight lattice of $T$ and 
let $\Lambda^+$\index{$\Lambda^+$, dominant integral weights} be the set of dominant integral weights.
For $\lambda\in \Lambda^+$ let $P_\lambda\supseteq B$\index{$P_\lambda$, parabolic subgroup} be a parabolic subgroup such that $\lambda$\index{$\lambda$, dominant integral weight} extends to a character of $P_\lambda$
and the associated line bundle $\mathcal L_\lambda$\index{ $\mathcal L_\lambda$, line bundle} on $G/P_\lambda$ is ample.
For a dominant weight $\lambda$ let $V(\lambda)\simeq H^0(G/P_\lambda,\mathcal L_\lambda)^*$ \index{$V(\lambda)$, irreducible $G$-representation} be the irreducible $G$-representation
of highest weight $\lambda$ and fix a highest weight vector $v_\lambda$\index{$v_\lambda$, highest weight vector}. 
Denote by $\iota:G/P_\lambda\hookrightarrow \mathbb P(V(\lambda))$ the embedding of the generalized flag variety $G/P_\lambda$\index{$G/P_\lambda$, generalized flag variety}
as the highest weight orbit $G.[v_\lambda]\subset \mathbb P(V(\lambda))$ and let $R=\mathbb C[G/P_\lambda]$ \index{$R$, homogeneous coordinate ring of $G/P_\lambda$} be the graded homogeneous
coordinate ring of the embedded variety. As a $G$-representation, the latter is isomorphic to
$\bigoplus_{n\ge 0} H^0(G/P_\lambda,\mathcal L_{n\lambda})$.
It is natural to ask: does Theorem~\ref{grassmanncase} have an appropriate reformulation in this vastly more general setting?

\subsection{The generalization and the price to pay} It turns out that Theorem~\ref{grassmanncase} holds (roughly) without any changes. The price one has to pay is that the notation becomes
heavier and the construction of the basis given by the standard monomials is not anymore as explicit as in the previous section:
the indexing set $I_{d,n}$ is replaced by the set of LS-paths of shape $\lambda$ (see Definition~\ref{lspath} below or \cite{L1}, LS-path is an abbreviation for 
Lakshmibai-Seshadri path),
the basis of $(\Lambda^d\mathbb C^n)^*$ given by the Pl\"ucker coordinates is replaced by the path vectors (see Definition~\ref{pathvector} below or \cite{L2}),
we still have relations called Pl\"ucker relations, and the
Hodge algebra theory is replaced by the theory of {\it LS-algebras}, a generalization of the Hodge algebras introduced by
Rocco Chiriv\`i \cite{Ch1}.

The description we give now is sometimes a little sloppy, we try more to explain where this construction comes from than
to dwell in precise technical details, these can be found in the corresponding articles.

\subsection{LS-paths}  Let $W$\index{$W$, Weyl group} be the Weyl group of $G$, it comes naturally endowed with a partial order "$>_B$"\index{$>_B$, Bruhat order} called the Bruhat order,
and a length function $\ell:W\rightarrow\mathbb N$. The value $\ell(w)$\index{$\ell(w)$, length function on $W$} can be defined as the dimension of the
associated Schubert variety $X(w)$ \index{$X(w)$, Schubert variety} in $G/B$, or as the minimal number of terms needed to write $w$ as a product of
simple reflections. Let $N$ \index{$N$, number of positive roots} be the number of positive roots of $G$.
A {\it maximal chain in $W$} is a linearly ordered sequence $\underline{\bf w}=(w_0,w_1,\ldots,w_N)$ \index{$\underline{\bf w}$, linearly ordered sequence of Weyl group elements} of Weyl group elements such that
$w_0>_Bw_1>_B\ldots>_Bw_N$ and $\ell(w_i)=N-i$. So $w_0$\index{$w_0$, longest Weyl group element} is always the unique element in $W$ of maximal length and $w_N$ is always the identity element.
Such a maximal chain gives rise to a sequence $(\beta_1,\ldots,\beta_N)$ \index{$\beta$, positive root} of positive roots such that $s_{\beta_j}w_{j-1}=w_j$.

\par\noindent
\begin{exam}\rm
For $G=SL_3$ the Weyl group $W$ is the symmetric group $\mathfrak S_{3}$ generated by the two simple reflections $s_1$ and $s_2$.
The Bruhat order can be recovered from the following four maximal chains:
$$
\begin{array}{ccc}
w_0>_Bs_1s_2 >_Bs_1>\text{id}, & w_0>_Bs_2s_1>_Bs_2>_B\text{id},\\
w_0>_Bs_1s_2>_Bs_2>\text{id}, &w_0>_Bs_2s_1 >_Bs_1>_B\text{id},
\end{array}
$$
leading to the following four sequences of roots:
$$
\begin{array}{ccc}
(\alpha_2, \alpha_1+\alpha_2, \alpha_1)&(\alpha_1, \alpha_1+\alpha_2, \alpha_2) \\
( \alpha_2, \alpha_1,  \alpha_2) & (\alpha_1, \alpha_2, \alpha_1).
\end{array}
$$
\end{exam}
Let $\lambda$ be a dominant weight and let $W_\lambda$ be the stabilizer of $\lambda$ in $W$. The
length function $\ell$, the Bruhat order, the notion of a maximal chain, the associated sequence of positive roots etc.,
all this can be defined for the quotient $W/W_\lambda$ too.
Let $N_\lambda$\index{$N_\lambda$, number of positive roots in $U_\lambda \subset P_\lambda$} be the number of positive roots corresponding to the unipotent radical of the parabolic subgroup $P_\lambda$.

\begin{definition}\label{LS-chain}\it
Let $\lambda$ be a dominant weight. An LS-chain of shape $\lambda$
is a pair $\hat\pi=(\underline{\bf w},\underline{\bf a})$ \index{$\hat\pi$, LS-chain}, where $\underline{\bf w}=(w_0,w_1,\ldots,w_{N_\lambda})$ is a maximal chain in $W$,
$(\beta_1,\ldots,\beta_{N_\lambda})$
is the associated sequence of positive roots and
$$
\underline{\bf a}=(a_{-1}=0\le a_0\le a_1\le\ldots\le a_{N_\lambda}=1)
$$
is a weakly increasing sequence of rational numbers such that
\begin{equation}\label{integrally}
\forall i=0,\ldots,N_\lambda-1: a_i\langle w_{i}(\lambda),\beta^\vee_i\rangle\in\mathbb Z. 
\end{equation}
\index{$\beta^\vee$, coroot}
\end{definition}
Let $\Lambda_{\mathbb R}=\Lambda\otimes_{\mathbb Z}\mathbb R$\index{$\Lambda_{\mathbb R}$, real span of the weight lattice} be the real span of the weight lattice.
Given an LS-chain $\hat\pi=(\underline{\bf w},\underline{\bf a})$,
then we define the associated path as follows, we set $a_{-1}=0$:
\begin{definition} \it
The associated {\it LS-path $\pi$}\index{$\pi$, LS-path} is the piecewise linear map:
$$
\pi:[0,1]\rightarrow \Lambda_{\mathbb R},\quad t\mapsto
\sum_{j=0}^{m-1} (a_i-a_{i-1}) w_{i}(\lambda) + (t-a_{m-1})w_{m}(\lambda)\quad\hbox{for\ } t\in[a_{m-1},a_m].
$$
\end{definition}
Different LS-chains may give rise to the same LS-path, up to reparametrization:
\par\noindent
\begin{exam}\rm
Let $G=SL_3$, fix $\lambda=\rho$ and set $\underline{\mathbf a}= (0\le 1\le 1\le 1\le 1)$.
The four LS-chains of type $\rho$: $((w_0,s_2s_1,s_1,\text{id}),\underline{\mathbf a})$, $((w_0,s_2s_1,s_2,\text{id}),\underline{\mathbf a})$,
$((w_0,s_1s_2,s_1,\text{id}),\underline{\mathbf a})$
and $((w_0,s_1s_2,s_2,\text{id}),\underline{\mathbf a})$ are LS-chains
giving rise to the same path: the straight line joining the origin with $w_0(\rho)$.
\end{exam}
Let $\hat\pi$ be an LS-chain, the
associated LS-path $\pi$ is obtained from $\hat\pi$ by omitting those entries for which the rational numbers are equal.
More precisely:
\begin{definition}\label{lspath} Let
$\hat\pi=(\underline{\bf w},\underline{\bf a})$ be an LS-chain. For the weakly increasing sequence
$\underline{\mathbf a}=(0\le a_0\le a_1\le\ldots\le a_{N_\lambda}=1)$ let
$0< i_1 < i_2<\ldots <i_k\le N_{\lambda}$ be the indices such that the value of the $a_i$ jumps, i.e.:
$$
0=a_{-1}=a_0=\ldots=a_{i_1-1}<a_{i_1}=\ldots=a_{i_2-1}<a_{i_2}=\ldots <a_{i_k}=\ldots=1.
$$
The LS-path $\pi$ of shape $\lambda$ associated to $\hat\pi$ is the pair of sequences
\begin{equation}\label{path}
\pi=((w_{i_1},\ldots,w_{i_k}),(0<a_{i_1}<a_{i_2}<\ldots< a_{i_k}=1)).
\end{equation}
\end{definition}
\par\noindent
\begin{exam}\rm
For the LS-chains $\hat\pi$ considered in Example 3 we get $0<a_0=a_1=a_2=a_3=1$ and hence
$\pi=((s_1s_2s_1),(0<1))$.
\end{exam}
It is now easy to see that Definition~\ref{lspath} is equivalent to the definition of an LS-path in \cite{L1}.
In fact, the definition given there means exactly that a pair of sequences for $\pi$ as above in \eqref{path} is an
LS-path if and only if it can be extended to an LS-chain in the sense of Definition~\ref{LS-chain}.

\begin{remark}
One should think of a maximal chain as a geometric object, this is also the way how Lakshmibai, Musili and Seshadri
(see, for example, \cite{LMS,LaSc1,LaSc2}) have been coming up with their combinatorics related to standard monomial theory,
which culminated in the conjectures about what now is called an LS-path. A maximal chain in $W$ corresponds to a maximal sequence of Schubert varieties in $G/P_\lambda$,
which are subsequently contained in each other. One should think of an LS-chain $\hat\pi$ as an object
corresponding to a section $s_\pi\in H^0(G/P_\lambda,\cal L_{\lambda})$ characterized by the ``vanishing behavior" on the Schubert varieties
occurring in this maximal chain. The vanishing multiplicities should be given by the integers $a_i\langle w_{i}(\lambda),\beta^\vee_i\rangle$.
The path vectors defined below mimic this behavior, up to filtration \cite{L2}.
For a better understanding of this interplay between algebraic combinatorics, geometry and representation theory,
 it would be very good to have an explicit bijection between the LS-chains/LS-paths and a basis of $H^0(G/P_\lambda,\cal L_{\lambda})$
with vanishing properties precisely described by the LS-paths, and which characterize the section.
\end{remark}
\subsection{Path-vectors} For a dominant weight $\lambda$, fix a lowest weight vector $p_{\text{id}}\in H^0(G/P_\lambda,\cal L_{\lambda})$.
For an element $w\in W/W_\lambda$ choose an appropriate lift $\hat w\in G$ and set $p_w:= \hat w(p_{\text{id}})$. Such a vector is uniquely defined up to a
scalar multiple and is called an extremal weight vector.
The following definition
does not make sense, but it describes quite well what kind of properties the path vector should have:
\begin{definition}\label{pathvector}\it
Let $\hat\pi=(\underline{\bf w},\underline{\bf a})$ be an LS-chain. The path vector $p_\pi$ associated to $\pi$ is the section
in $H^0(G/P_\lambda,\cal L_{\lambda})$ defined by:
$$
p_\pi:={}^q\sqrt{p_{w_0}^{qa_0} p_{w_1}^{q(a_1-a_0)}\cdots p_{w_{N_\lambda}}^{q(1-a_{N_\lambda-1})}}
$$
\index{$p_\pi$, path vector} where $q\in \mathbb Z_{>0}$ is minimal such that $q a_i\in\mathbb Z$ for all $i$.
\end{definition}
{\bf Some explanation:} In the definition we switch from the LS-chain $\hat \pi$ to the LS-path $\pi$. It is easy to see that the product
$p_{w_0}^{qa_1} p_{w_1}^{q(a_2-a_1)}\cdots p_{w_{N_\lambda}}^{q(1-a_{N_\lambda})}$ of extremal weight vectors is independent
of the chain we started with, \emph{i.e.} if two chains lead to the same path, then the corresponding products
are the same.

It remains to say something about ``${}^q\sqrt{\ }$''.
The product of the extremal weight vectors makes sense,
this is a section living in $H^0(G/P_\lambda,{\cal L}_{q\lambda})$. It remains to explain how to ``take its $q$-th root", i.e., how
to produce out of the given section one in $H^0(G/P_\lambda,\cal L_{\lambda})$.
This is described in detail in \cite{L2}, where a quantized version of the Frobenius splitting is constructed,
providing a map from $H^0(G/P_\lambda,{\cal L}_{q\lambda})_\xi$ (i.e. the module is viewed as a representation for the quantum group at a
$q$-th root of unity $\xi$) to $H^0(G/P_\lambda,{\cal L}_{\lambda})$
see also \cite{KuL}. This splitting can be viewed as a procedure to take a $q$-th root out of a section.

\subsection{Standard monomials}
Let $\lambda\in \Lambda^+$ be a dominant weight and let $\pi=((w_{i_1},\ldots,w_{i_k}),(0<a_{i_1}<a_{i_2}<\ldots< a_{i_k}=1))$
be an LS-path of shape $\lambda$. We set $i(\pi)=w_{i_1}$ (the initial direction of the path) and
 $e(\pi)=w_{i_k}$ (the final direction of the path).
\begin{definition}\label{standardmonomial2}
Let $\pi_1,\ldots,\pi_s$ be LS-paths of shape $\lambda$. The product $p_{\pi_1}\cdots p_{\pi_s}$
of the corresponding section is called a {\it standard monomial} in $H^0(G/P_\lambda,{\cal L}_{s\lambda})$
if  $e(\pi_1)\ge i(\pi_2)\ge e(\pi_2)\ge\ldots\ge i(\pi_s)$.
\end{definition}
We say that two LS-paths $\pi_1,\pi_2$ of shape $p\lambda$ respectively $q\lambda$ have the same support
if there exist two LS-chains $\hat \pi_1=(\underline{\bf u},\underline{\bf a}),\pi_2=(\underline{\bf v},\underline{\bf b})$
(corresponding to $\pi_1$ respectively $\pi_2$) of shape $p\lambda$ respectively $q\lambda$ such that
$\underline{\bf u}=\underline{\bf v}$. If we have such a pair, then we construct a new LS-chain $\hat \pi =(\underline{\bf u},\underline{\bf c})$
of shape $(p+q)\lambda$, where the $c_j$ are defined inductively as follows: $c_{-1}=0$, and for $j\ge 0$ we set
$c_j=\frac{p(a_j-a_{j-1})+q(b_j-b_{j-1})}{p+q}$. Note that $\hat\pi$ is an LS-chain, denote by $\pi$ the associated LS-path.
One can define on the set of all LS-paths a monomial order (a weighted lexicographic order).
For the following see  \cite{L2}:
\begin{theorem}
\begin{enumerate}
\item[{\it i)}] The standard monomials form a basis of  the ring
$$R=\bigoplus_{n\ge 0} H^0(G/Q,\mathcal L_{n\lambda}).$$
\item[{\it ii)}] Special Pl\"ucker relations:  Let $\pi_1,\pi_2$ be LS-paths of shape $p\lambda$ respectively $q\lambda$ having the same support.
Then $p_{\pi_1}p_{\pi_2}=p_\pi + \sum_{\eta>\pi} a_{\eta} p_\eta$, where $\pi$ is constructed out of $\pi_1,\pi_2$ as above.
\end{enumerate}
\end{theorem}
This theorem, which can be viewed as a generalization of Theorem~\ref{grassmanncase}, was the starting point
for Chiriv\`i  \cite{Ch1} to generalize the notion of an algebra with straightening laws.
Applying this to the embedded flag variety, his theory implies:
\begin{theorem}
There exists a flat degeneration of $G/P_\lambda\hookrightarrow \mathbb P(H^0(G/P_\lambda,\mathcal L_{\lambda})^*)$
into a union of toric varieties $\cal T$\index{$\mathcal{T}$, union of toric varieties}, one irreducible component for each maximal chain in $W/W_\lambda$.
The homogeneous coordinate ring of $\cal T$ has a basis $\{p_\eta\}$ indexed by all LS-paths of shape
$n\lambda$, $n\in\mathbb N$, with multiplication rule: $p_{\pi_1}p_{\pi_2}=0$ if they have no common support,
and $p_{\pi_1}p_{\pi_2}=p_\pi$ if they have a common support, and $\pi$ is constructed out of $\pi_1,\pi_2$ as above.
\end{theorem}
\begin{remark}
For the extension to the multi-cone over $G/P_\lambda$ see \cite{Ch1},
for the degeneration of Schubert varieties (the standard monomial basis is compatible with
all Schubert varieties) see  \cite{Ch1,Ch2}.
\end{remark}
\subsection{Other degenerations using standard monomial theory}
\begin{example}\label{GonLak}
{\bf The Gonciulea-Lakshmibai degeneration.} \rm
It is natural to ask for flat toric degenerations such that the special fiber
remains irreducible.  The first degeneration of the flag variety $SL_n/B$ with this property was obtained by N. Gonciulea and V. Lakshmibai
in \cite{GL}, they have been using standard monomial theory. The essential point
in their proof is the fact that the fundamental weights
are minuscule weights. This can be used to show that the indexing system for a basis of every
fundamental representation is endowed with a structure of distributive lattice. Roughly speaking, this allows them to
degenerate the Pl\"ucker relation in \eqref{pluecker} to the equation
$\mathbf{p}_{\underline{\mathbf i}}\mathbf{p}_{\underline{\mathbf j}}=\mathbf{p}_{\underline{\mathbf i}\cup \underline{\mathbf j}} \mathbf{p}_{\underline{\mathbf i}\cap \underline{\mathbf j}}$
and not just to $\mathbf{p}_{\underline{\mathbf i}}\mathbf{p}_{\underline{\mathbf j}}=0$ for $\underline{\mathbf i},\underline{\mathbf j}$ not comparable,
as it is done using the theory of Hodge algebras.
For more details see \cite{GL}.
\end{example}
\begin{example} {\bf Generalization to Schubert varieties by Dehy and Yu.} \rm
A toric degeneration for Schubert varieties in the $SL_n$-case is described in \cite{DY2} by R. Dehy and R. W. T. Yu.
They use the polytopes described in \cite{DY1} and methods similar to those in \cite{GL}.
Using LS-paths for semisimple or affine Lie algebras of rank 2 (i.e. of types $\tt A_2$, $\tt B_2$, $\tt G_2$, $\tt A^{(1)}_1$ or $\tt A^{(2)}_2$),
R. Dehy \cite{De} proves that there exists a flat deformation for all Schubert varieties into a toric variety.
The methods are similar to those
used in \cite{DY2}.
\end{example}
\subsection{Irreducibility versus compatibility}
It is easy to see that Chiriv\`i's degeneration is an example for a degeneration which is compatible with all
Schubert varieties, their unions and intersections simultaneously. One may ask whether there exists a flat degeneration which is
also compatible with all Schubert varieties and, in addition, has an irreducible toric variety as a special fiber. Such a degeneration
does not exist, as was already pointed out by Olivier Mathieu: intersections of irreducible toric varieties
are irreducible toric varieties, but an intersection of Schubert varieties can be a union of
several Schubert varieties.

\section{The \texorpdfstring{variety ${G/\hskip -4.0pt/U}$}{flag variety} and Caldero's degeneration}\label{GnachU}
Again we start with a basis, this time it is the dual canonical basis of $\mathbb C[{G/\hskip -4.0pt/U}]$ arising from
the quantized enveloping algebra \cite{Ka,Lu1,Lu2}. Having in mind
the result by Chiriv\`i and the construction by Gonciulea and Lakshmibai for $G = SL_n$, one may ask:
Is there a flat family over $\hbox{Spec\,}\mathbb C[t]$ such that the generic fiber
is $G/P_\lambda$ and the special fiber is a toric variety (so in particular, irreducible)?
For $G=SL_n$ the answer is yes by \cite{GL}, see Example~\ref{GonLak}. A positive answer in the general case
was given by Caldero \cite{Ca1}, he uses special properties of the dual canonical basis
of $\mathbb C[G/\hskip -3.5pt/U]$.

\subsection{Some notation} Let $G$ be a connected complex reductive algebraic group
isomorphic to $G^{ss}\times (\mathbb C^*)^r$, where $G^{ss}$\index{$G^{ss}$, semisimple part of $G$} denotes the semisimple part of $G$,
and $G^{ss}$ is simply connected. Let $\hbox{Lie\,}G=\mathfrak g$ \index{$\mathfrak{g}$, Lie algebra } be its Lie algebra.
We fix a Cartan decomposition $\mathfrak g=\mathfrak n^-\oplus\mathfrak b$ \index{$\mathfrak n^-$, opposite nilpotent radical} \index{$\mathfrak b$, Borel subalgebra}, where
$\mathfrak b$ is a Borel subalgebra with maximal torus $\mathfrak t$\index{$\mathfrak t$, maximal torus of $\mathfrak{b}$} and nilpotent radical
$\mathfrak n^+$\index{$\mathfrak n^+$, nilpotent radical}.  Let $\Phi$\index{$\Phi$, root system} be the root system of $\mathfrak g$ and denote by $\Phi^+$\index{$\Phi^+$, positive roots} the set of
positive roots with respect to the choice of $\mathfrak b$. We denote by ${\cal R}$\index{ $\cal R$, root lattice} the root lattice and by
${\cal R}^+$\index{${\cal R}^+$, semi-group of positive roots} the semigroup generated by the positive roots. For $\beta\in \Phi^+$, let
$f_\beta$\index{$f_\beta$, root vector} be a non-zero root vector in $\mathfrak n^-$ of weight $-\beta$.

Let $B\subset G$ \index{$B$, Borel subgroup} be the Borel subgroup such that $\hbox{Lie\,}B=\mathfrak b$, let $U\subset B$\index{$U$, unipotent radical} be its unipotent radical
and denote by $T\subset B$ the maximal torus such that $\hbox{Lie\,}T=\mathfrak t$.
We write ``$\succ_{wt}$''\index{$\succ_{wt}$, partial order on the weight lattice} for the usual partial order on $\Lambda$, i.e., $\lambda\succ_{wt}\mu$ if and only if $\lambda-\mu\in{\cal R}^+$.
We denote by $B^-$ the opposite Borel subgroup and let $U^-$ \index{$U^-$, opposite unipotent radical} be its unipotent radical, so $\hbox{Lie\,}U^-=\mathfrak n^-$.
Let $n$\index{$n$, rank of $\Phi$} be the rank of $\Phi$ and $n+r$\index{$n+r$, rank of $G$} be the rank of $G$.

\subsection{The \texorpdfstring{variety ${G/\hskip -4pt/U}$}{flag variety}}\label{flagvariety}
 One has a natural $G\times G$-action on $G$ by left and right multiplication.
The ring $\mathbb C[G]^{1\times U}$\index{ $\mathbb C[G]^{U}$. $1 \times U$ invariant functions} of $1\times U$-invariant functions (in the following we write just $\mathbb C[G]^{U}$) is
finitely generated and normal, so the variety $G/\hskip -3.5pt/U:=\hbox{\textrm{Spec}}(\mathbb C[G]^U)$\index{ $G/\hskip -3.5pt/U:=\hbox{\textrm{Spec}}(\mathbb C[G]^U)$, $U$-invariant points}  is a normal
(but in general singular) affine variety.
Since $1\times T$ normalizes $1\times U$, $\mathbb C[G]^U$ is a natural $G\times T$-algebra.
As a $G$-representation (resp. $G\times T$-representation), its coordinate ring is isomorphic to
\begin{equation}\label{coordinateGnachU}
\mathbb C[G/\hskip -3.5pt/U]\simeq \bigoplus_{\lambda\in\Lambda^+} V(\lambda)^*\simeq \bigoplus_{\lambda\in\Lambda^+} V(\lambda)^*\otimes v_\lambda\simeq \mathbb C[G]^U.
\end{equation}
So  $(1,t)\in G\times T$ acts on $V(\lambda)^*\simeq V(\lambda)^*\otimes v_\lambda$ by the scalar $\lambda(t)$.

The variety $G/\hskip -3.5pt/U$ is endowed with a natural $G$-action, making it into a spherical variety, i.e., a variety with a dense $B$-orbit.
One has a canonical dominant map $\psi:G\rightarrow G/\hskip -3.5pt/U$, inducing an inclusion $G/U\hookrightarrow G/\hskip -3.5pt/U$
and a birational orbit map  $o:B^-\rightarrow G/\hskip -3.5pt/U$, $b\mapsto b.\bar{1}$ (where $\bar{1}=\psi(1)$).

An element $\phi\in V(\lambda)^*$ can be seen as a function on the open and dense subset $G/U\hookrightarrow G/\hskip -3.5pt/U$
as follows: for a class $\bar g\in G/U$ let $g\in G$ be a representative. Then \index{$\phi\vert_{G/U}$, coordinate function}
\begin{equation}\label{coordinatefunction}
\phi\vert_{G/U}: G/U\rightarrow \mathbb C,\quad \bar g\mapsto \phi(g v_\lambda).
\end{equation}

\subsection{The string cone}\label{thestringcone}
The algebra $\mathbb C[G/\hskip -3.5pt/U]$ has a basis with some remarkable properties.
Recall that Kashiwara and Lusztig constructed a
global crystal basis respectively canonical basis (two names for the same basis, see \cite{GrL}) for representations $V_q(\lambda)$
of the corresponding quantum groups.
Specializing at $q=1$, the dual canonical basis $\mathbb B^*$ \index{$\mathbb B^*$, dual canonical basis} is the basis of
$\mathbb C[G/\hskip -3.5pt/U]=\bigoplus_{\lambda\in\Lambda^+}V(\lambda)^*$ dual to the canonical basis $\mathbb{B}$\index{$\mathbb{B}$, canonical basis}
of $\bigoplus_{\lambda\in\Lambda^+}V(\lambda)$.

The elements of the dual canonical basis are indexed by two parameters: a dominant weight
$\lambda$ and an $N$-tuple of integers, where the latter depends on
the choice of a reduced decomposition $\underline{w}_0=s_{i_1} s_{i_2}\cdots s_{i_N}$ \index{$\underline{w_0}$, reduced expression of $w_0$}
of the longest word $w_0$ in $W$. Let ${\cal B}\subset U_q(\mathfrak n^-)$ \index{ ${\cal B}$, canonical basis of $U_q(\mathfrak n^-)$} \index{$U_q(\mathfrak n^-)$, lower part of the quantum group}  be the canonical basis and denote
by ${\cal B^*}\subset U_q(\mathfrak n)$ \index{${\cal B^*}$, dual canonical basis} \index{$U_q(\mathfrak n)$, positive part of the quantum group} the dual canonical basis.
Using either Kashiwara operators on the crystal basis (see for example in \cite{L3}) or using quantum derivations of $U_q(\mathfrak n)$ (see for example \cite{Ca1}),
one defines a map ${\cal I}_{\underline{w}_0}: {\cal B^*}\rightarrow \mathbb N^N$ \index{${\cal I}_{\underline{w}_0}$, string coordinates map}, which associates to an element  $b^*\in {\cal B^*}$ \index{$b^*$, dual canonical basis element} its {\it string coordinates}.

It has been proved in \cite{L3}, see also \cite{BZ}, that the map is injective, and the image of ${\cal I}_{\underline{w}_0}$ is the monoid of integral
points of a rational convex polyhedral cone $C_{\underline{w}_0} \subset\mathbb R^N$, \index{$C_{\underline{w}_0}$, string cone} called the {\it string cone}.
Let now ${\cal B}_\lambda\subset\cal B$ \index{${\cal B}_\lambda$, canonical basis for $V_q(\lambda)$} be the subset of elements such that $b.v_\lambda\not=0$ in $V_q(\lambda)$.
One can identify ${\cal B}_\lambda$ with the canonical basis $\{b.v_\lambda\mid b\in {\mathbb B}_\lambda\}$
of $V_q(\lambda)$. Let $\hat{\cal B}$ be the disjoint union
of the ${\cal B}_\lambda$. In this way we can view $\hat{\cal B^*}$ as the disjoint union of the dual bases
$\mathbb B^*_\lambda$, and the parametrization of the basis elements is given by the rule that $b_{\lambda,\mathbf m}$\index{$b_{\lambda,\mathbf m}$, dual canonical basis element} is the unique
element in $\mathbb B^*_\lambda\subset {\cal B^*}$ such that ${\cal I}_{\underline{w}_0}(b_{\lambda,\mathbf m})=\mathbf m$.
Let
\begin{equation}\label{semigroupfirst}
{\Gamma}_{{\underline{w}_0}}=\{(\lambda,\mathbf m)\in\Lambda^+\times \mathbb N^N\mid \exists \lambda\in\Lambda^+,\,\exists
b\,\in \mathbb B^*_\lambda: {\cal I}_{\underline{w}_0}(b)=\mathbf m\}.
\end{equation}
\index{$\Gamma_{{\underline{w}_0}}$, semi-group for the string cone} It has been shown in \cite{L3}, see also \cite{BZ}, that ${\Gamma}_{{\underline{w}_0}}$ is the monoid of integral point
of a rational convex polyhedral cone ${\mathcal C}_{{\underline{w}_0}}$\index{${\mathcal C}_{{\underline{w}_0}}$, polyhedral string cone}.

So by the construction above, the algebra $\mathbb C[G/\hskip -3.5pt/U]$ comes equipped with the basis
$\mathbb B^*=\{b_{\lambda,\mathbf m}\mid (\lambda,\mathbf m)\in {\Gamma}_{{\underline{w}_0}}\}$. The basis
elements are $T\times T$-eigenvectors, where $b_{\lambda,\mathbf m}$ is of weight $\lambda$ for
the right action and of weight $-\lambda + m_1\alpha_{i_1}+\ldots + m_N\alpha_{i_N}$ for the left action of $T$.

\subsection{The degeneration}\label{thedegeneration}
A remarkable property of the dual canonical basis is the following multiplication rule
proved by Caldero:
$$
b_{\lambda,\mathbf m} b_{\mu,\mathbf n}= b_{\lambda+\mu,\mathbf{m+n}} +
\sum_{\mathbf{k}>\mathbf{m+n}}  c^\mathbf{k}_{(\lambda,\mathbf{m}),(\mu,\mathbf{n})} b_{\lambda+\mu,\mathbf{k}},
$$
where ``$\le$'' denotes the lexicographic ordering on $\mathbb N^N$. So the multiplication rule for basis elements
can be described as: up to elements which are larger with respect to the lexicographic ordering,
the product of basis elements is the same as in the monoid ${\Gamma}_{{\underline{w}_0}}$, which is the index system
of $\mathbb B^*$.

From this, Caldero deduces the existence of an increasing filtration of $\mathbb C[G/\hskip -3.5pt/U]$
by $T\times T$-submodules, such that the associated graded algebra, $gr\,\mathbb C[G/\hskip -3.5pt/U]$,
is isomorphic to the algebra of the monoid $\mathbb C[{\Gamma}_{{\underline{w}_0}}]$.
In geometric terms:
\begin{theorem}[\rm \cite{Ca1}]
The affine variety $G/\hskip -3.5pt/U$ admits a flat degeneration
to a normal affine toric variety $X_0={\rm Spec\,}\mathbb C[{\Gamma}_{{\underline{w}_0}}]$\index{$X_0$, special fiber} for the torus $T\times \mathbb T$,
where we put $\mathbb T := (\mathbb C^*)^N$.
Further, the degeneration is compatible with the actions of $T\times T$ on $G/\hskip -3.5pt/U$
(regarding $T\times T$ as a subgroup of $G\times T$), and on $X_0$ via the homomorphism
of tori $T\times T\rightarrow T\times\mathbb  T$, $(t, t')\mapsto (t^{-1}t', \alpha_{i_1} (t), . . . , \alpha_{i_N} (t))$.
\end{theorem}

\subsection{Degeneration of \texorpdfstring{${G/P_\lambda}$}{the flag variety} by Alexeev and Brion}\label{ALBriondeg}
For a given dominant weight $\lambda$ and a fixed reduced decomposition $\underline{w}_0$,
the set of integral points
$$
\{{(\mathbf m})\mid \exists b\,\in {\cal B}^*_\lambda: {\cal I}_{\underline{w}_0}(b)=\mathbf m\}\subset \mathbb R^N
$$
is in bijection with the set of integral points of a rational convex polytope $Q_{\underline{w}_0}(\lambda)$ \index{$Q_{\underline{w}_0}(\lambda)$, string polytope},
called the {\it string polytope} for $\lambda$ \cite{BZ,L3}. This polytope has a geometric interpretation. The following
theorem holds (appropriately reformulated) in much more generality for polarized spherical varieties and not only
for flag varieties.
\begin{theorem}{\rm \cite{AB}}
There exists a family of $T$-varieties $\pi: {\cal X} \rightarrow \mathbb A^1$, where ${\cal X}$ is a
normal variety, such that $\pi$ is projective and flat, it is trivial with fiber $G/P_\lambda$ over the complement of $0$ in $\mathbb A^1$,
the fiber of $\pi$ at $0$ is isomorphic to $X_0$, which is a toric variety for the torus $\mathbb T$. The string polytope
$Q_{\underline{w}_0}(\lambda)$ can be recovered as the moment polytope of the
toric variety $X_0$.
\end{theorem}
\subsection{Kaveh's interpretation of \texorpdfstring{$Q_{\underline{w}_0}(\lambda)$}{the string polytopes} as Newton-Okounkov body} For more
details on Newton-Okounkov bodies we refer to \cite{KK}, see also section~\ref{NObody}. 
Kaveh proves in \cite{K1} that the string
parametrization of the crystal basis mentioned at the beginning of this section
coincides with a natural valuation on the field of rational functions on the flag variety $G/B$.
The $\mathbb Z^N$-valued valuation can be defined using a coordinate system on a Bott-Samelson variety,
so it depends on the choice of a reduced decomposition ${\underline{w}_0}$ of the longest word $w_0$
in the Weyl group. It turns out that the associated monoid (see also section~\ref{Section:birational}) is exactly ${\Gamma}_{{\underline{w}_0}}$,
hence giving ${\Gamma}_{{\underline{w}_0}}$ an additional algebraic geometric interpretation.
This shows that the string polytopes associated to irreducible representations, can be realized as
Newton-Okounkov bodies for the flag variety. This fully generalizes an earlier result of
A. Okounkov \cite{O} for the generalized Gelfand-Tsetlin polytopes of the symplectic group.

\subsection{Anderson's connection between N-O-bodies and toric degeneration}\label{AndersonNO}
We have seen above that the string polytope $Q_{\underline{w}_0}(\lambda)$ can be viewed as a 
Newton-Okounkov body for the generalized flag variety $G/P_\lambda$ but also as the moment 
polytope for a toric variety $X_0$ obtained by a flat degeneration from $G/P_\lambda$. Dave Anderson \cite{A} has shown 
that this is always true once knows that the monoid used in the Newton-Okounkov approach is finitely generated.

\section{Newton-Okounkov bodies for flag varieties}\label{Section:birational}

\subsection{Introduction} We present a different approach to construct $T$-equi\-variant flat toric 
degenerations of flag varieties (see \cite{FaFL}). The procedure uses representation theory as well as ideas 
from the Newton-Okounkov theory. The strategy can be seen as a common generalization of Caldero's degeneration 
(and the subsequent constructions of flat toric degenerations of flag varieties by Alexeev-Brion \cite{AB}  
and Kaveh \cite{K1}) and the construction of flat toric degenerations using a refinement
of the PBW-filtration \cite{FFL1} (see also \cite{FFL2,FFL3}). We will see that both:
Caldero's construction and the construction via filtrations have implicitly in common the idea to replace the group 
$U^-$ by a birationally equivalent product of root subgroups.   

To be more precise: to construct monoids like $\Gamma_{\underline{w}_0}$ (see \eqref{semigroupfirst} in section~\ref{thestringcone})
and associated toric degenerations of $G/\hskip -3.5pt/U$ with tools coming from the theory of Newton-Okounkov bodies, we need 
$\mathbb Z^{N+n}$-valued valuations on the field of rational functions $\mathbb C(G/\hskip -3.5pt/U)$.
Here $N$ is the number of positive roots in the root system $\Phi$ of $G$ and $n$ is the rank of $G$.
To construct such valuations, we consider sequences $(\beta_1,\ldots,\beta_N)$ of positive roots
(not necessarily pairwise different), such that the product map 
$$
U_{-\beta_1}\times\cdots\times U_{-\beta_N}\rightarrow U^-
$$ 
is birational (see section~\ref{birat} for examples). By fixing a basis of the weight lattice of $T$,
the birational maps $B^-\rightarrow G/\hskip -3.5pt/U$ and 
$$
T\times U_{-\beta_1}\times\cdots\times U_{-\beta_N}\rightarrow B^-=TU^-
$$
provide a natural identification of $\mathbb C(G/\hskip -3.5pt/U)$ with a function field in $N+n$ variables.
After fixing a total order on $\mathbb Z^{N+n}$ (see section~\ref{lexmonorder} for examples), the lowest term valuation 
(see section~\ref{SValuationsemigroup}) gives 
us then a natural $\mathbb Z^{N+n}$-valued valuation $\nu$ on $\mathbb C(G/\hskip -3.5pt/U)$. 
We associate to $\mathbb C[G/\hskip -3.5pt/U]$ the {\it valuation monoid} (see \eqref{valuationsemigroup}):
$$
{\cal V}(G/\hskip -3.5pt/U)=\{\nu({p})\mid p\in \mathbb C[G/\hskip -3.5pt/U]-\{0\}\},
$$ 
which depends of course on the choice of the sequence of positive roots and the choice of the total order.

This sequence of positive roots and the fixed total order can also be used to define a filtration
of the enveloping algebra $U(\mathfrak n^-)$ and on the finite dimensional irreducible highest weight 
representations of $G$. We use these filtrations to introduce the term of an {\it essential monomial}
for a representation and we associate to 
$G/\hskip -3.5pt/U$ the {\it global essential monoid} (see section~\ref{Filtessent})
$$
\Gamma= \{(\lambda,\mathbf m)\mid \lambda \hbox{\it\ a dominant weight},\ \mathbf f^{\mathbf m}
\hbox{\it \,essential for $V(\lambda)$}\} \subset \Lambda^+\times \mathbb N^N.
$$
This approach is close to the one which has been used in \cite{F1,FFL2,FFL3,FFL1}. It turns out that for a 
fixed sequence of positive roots and fixed choice of the total order (see section~\ref{gammagleichnu}): 
$$
\Gamma={\cal V}(G/\hskip -3.5pt/U).
$$
The proof shows that once one knows that the monoid $\Gamma$ is finitely generated and saturated,
then the methods discussed in sections~\ref{thedegeneration}--\ref{AndersonNO} directly apply also 
to this much more general situation.  
In particular, the degenerations of flag varieties by Caldero \cite{Ca1} and Alexeev and Brion \cite{AB},
as well as the degenerations via filtrations \cite{FFL1} show up in a uniform framework.

For examples and explicit descriptions of associated Newton-Okounkov bodies, 
questions about unimodular (in-) equivalence of the bodies etc., see section~\ref{examplesectionASM}. 
For an application in symplectic geometry see section~\ref{Gromov}.

The construction can be extended to spherical varieties in the same way as in \cite{AB} and \cite{K1}. 

\subsection{Birational sequences and examples}\label{birat}
As an affine variety, the subgroup $U^-\subset G$ is isomorphic
to $\mathbb A^N$. For a root vector $f_{\beta}\in {\mathfrak g}_{-\beta}$,
$\beta\in\Phi^+$, let $U_{-\beta}$
be the corresponding root subgroup $U_{-\beta}:=\{\exp(s f_{\beta})\mid s\in \mathbb C\}$ of $G$.

Fix a sequence of positive roots $S=({\beta_1},\ldots,{\beta_N})$.
We make no special assumption on this sequence, for example there may be repetitions, see
Example~\ref{monomialszwo}. Let $\mathbb T$ be the torus $(\mathbb C^*)^N$,
we write $\texttt{t}=(\texttt{t}_1,\ldots,\texttt{t}_N)$ for an element of $\mathbb T$.
\begin{definition} \label{orderedgenerating}
The variety ${Z}_S$ is the affine space $\mathbb A^N$ endowed with the following $T\times\mathbb T$-action:
$$
\forall (t,\texttt{t})\in T\times \mathbb T:(t,\texttt{t})\cdot(z_1,\ldots,z_N):=
(\texttt{t}_1 \beta_1(t)^{-1}z_1,\ldots, \texttt{t}_N \beta_N(t)^{-1}z_N).
$$
We call $S$\index{$S$, birational sequence}  a {\it birational sequence} for $U^-$ if the product map $\pi$ is birational:
\begin{equation}\label{birational1}
\pi:Z_S\rightarrow U^-,\quad (z_1,\ldots,z_N)\mapsto \exp(z_1f_{\beta_1})\cdots \exp(z_Nf_{\beta_N}).
\end{equation}
\end{definition}

We provide some examples of birational sequences, details and more examples can be found in \cite{FaFL}.

\begin{exam}\label{monomialseins} {\it The PBW-type case\/}: \rm
Let $S=({\beta_1},\beta_2,\ldots,{\beta_N})$ be an enumeration of the set of positive roots,
i.e. $\Phi^+=\{\beta_1,\ldots,\beta_N\}$. Then $S$ is a birational sequence.
\end{exam}
\begin{exam}
\label{monomialszwo} {\it The reduced decomposition case\/}: \rm
Fix a reduced decomposition $\underline{w}_0=s_{i_1}\cdots s_{i_N}$ of the longest word in the Weyl group and
a sequence $S=({\alpha_{i_1}},\ldots, {\alpha_{i_N}})$. Then $S$ is a birational sequence.
\end{exam}


\begin{example}\rm
If $S=({\beta_1},\ldots,{\beta_N})$ is a birational sequence, then so is $S'=({\beta_N},\ldots,{\beta_1})$.
\end{example}
\begin{example}\rm
Here is a list of all birational sequences for $G=SL_3$: $(\alpha_1,\alpha_2,\alpha_1)$, $(\alpha_2,\alpha_1,\alpha_2)$, $(\alpha_1,\alpha_1+\alpha_2,\alpha_2)$, $(\alpha_2,\alpha_1+\alpha_2,\alpha_1)$, $(\alpha_1,\alpha_2,\alpha_1+\alpha_2)$, $(\alpha_2,\alpha_1,\alpha_1+\alpha_2)$, $(\alpha_1+\alpha_2,\alpha_2,\alpha_1)$, $(\alpha_1+\alpha_2,\alpha_1,\alpha_2)$.
\end{example}

Let $S=({\beta_1},\ldots,{\beta_N})$ be a birational sequence for $U^-$.
Let ${\cal Z}_S$ be the toric variety ${\cal Z}_S= Z_S\times T$, where the torus
$T\times \mathbb T$ is acting on ${\cal Z}_S$ as follows:
$$
\forall (t,\texttt{t})\in T\times \mathbb T:(t,\texttt{t})\cdot(z,t'):=
(\texttt{t}_1 \beta_1(t)^{-1}z_1,\ldots, \texttt{t}_N \beta_N(t)^{-1}z_N;tt').
$$
Let $x_i:\, (z_1 ,\ldots,z_N)\mapsto z_i$ be the $i$-th coordinate function on ${Z}_S$,
then
$$
\mathbb C[{\cal Z}_S]\simeq \mathbb C[x_1,\ldots,x_N]\otimes \mathbb C[T]=\mathbb C[x_1,\ldots,x_N]\otimes\mathbb C[e^{\lambda}\mid\lambda\in\Lambda].
$$
Moreover, the canonical map ${\cal Z}_S\rightarrow G/ U$, $(z,t)\mapsto
\pi(z)\cdot t\cdot \bar{1}$, induces a $T$-equivariant birational map $\varphi: {\cal Z}_S\rightarrow G/\hskip -3.5pt / U$.
Via the isomorphism $\mathbb C(G/\hskip -3.5pt/U)\simeq\mathbb C({\cal Z}_S)$,
we can identify the coordinate ring $\mathbb C[G/\hskip -3.5pt/U]$ with a subalgebra of
$A:=\mathbb C[x_1,\ldots,x_N]\otimes\mathbb C[e^{\lambda}\mid\lambda\in\Lambda^+]\subset \mathbb C({\cal Z}_S)$.

\subsection{Weight functions and lexicographic orders}\label{lexmonorder} Let $\Psi: \mathbb Z^N\rightarrow \mathbb Z$
be a $\mathbb Z$-linear map such that $\Psi(\mathbb N^N)\subseteq \mathbb N$, we call $\Psi$\index{$\Psi$, integral weight function}
an  {\it integral weight function}. The {\it weight order} on $\mathbb N^N$ associated to $\Psi$ is the partial
order defined by $\mathbf m>_{\Psi} \mathbf m'$ iff $\Psi(\mathbf m)>\Psi(\mathbf m')$.
We write $>_{lex}$ for the lexicographic order and $>_{rlex}$
for the right lexicographic order on $\mathbb N^N$.
We refine the partial order above to a total order by combining the partial weight order
and the lexicographic order. In the following we use only
the lexicographic order, the generalization to the right lexicographic order is obvious.
\begin{definition}\label{weightorder}
A {\it $\Psi$-weighted lexicographic order} on
$\mathbb N^N$ is a total order ``$>$" on $\mathbb N^N$ refining ``$>_{\Psi}$'' as follows:
$$
\mathbf m>\mathbf m'\Leftrightarrow \hbox{\ either $\Psi(\mathbf m)> \Psi(\mathbf m')$, or
$\Psi(\mathbf m)=\Psi(\mathbf m')$ and $\mathbf m>_{lex}\mathbf m'$ }.
$$
If $\Psi$ satisfies in addition the condition $\Psi(\mathbf m)>0$ for all $\mathbf m\in \mathbb N^N-\{0\}$,
then a {\it $\Psi$-weighted opposite lexicographic order} on
$\mathbb N^N$ is a total order ``$>$" on $\mathbb N^N$ defined as follows:
$$
\mathbf m>\mathbf m'\Leftrightarrow \hbox{\ either $\Psi(\mathbf m)> \Psi(\mathbf m')$, or
$\Psi(\mathbf m)=\Psi(\mathbf m')$ and $\mathbf m<_{lex}\mathbf m'$}.
$$
\end{definition}
It is easy to verify that a $\Psi$-weighted lex order defines a monomial order on $\mathbb N^N$, \emph{i.e.}, a total order such that $\bf m>m'$ implies for all ${\bf m''}\in \mathbb N^N-\{0\}$: $\bf m+m''>m'+m''>m'$.

\begin{exam}\label{homorder} \rm
If $\Psi:\mathbb Z^N\rightarrow \mathbb Z$ is the map $\mathbf m\mapsto \sum_{i=1}^N m_i$, then the $\Psi$-weighted lexicographic order
is the homogeneous lexicographic order.
\end{exam}
\begin{exam}\label{rootorder} \rm
Fix a sequence $S=(\beta_1,\ldots,\beta_N)$ of roots in $\Phi^+$, let $ht$ be the height function on the positive roots and let
$\Psi$ be the {\it height weighted function}:
$$
\Psi:\mathbb Z^N\rightarrow \mathbb Z,\quad \mathbf m\mapsto \sum_{i=1,\ldots,N} m_i ht(\beta_i).
$$
\end{exam}
\subsection{The valuation monoid}\label{SValuationsemigroup}
Let $\omega_1,\ldots, \omega_{n'}$ be the fundamental weights of the semisimple part of $G$ and extend
this set by $\eta_{1},\ldots, \eta_{n-n'}$ to a basis of the character group of $T$.
Set $\Lambda^\dagger=\mathbb N\omega_1\oplus\ldots \oplus \mathbb N\omega_{n'}\oplus
\mathbb N\eta_{1}\oplus\ldots \oplus \mathbb N\eta_{n-n'}\subset \Lambda^+$\index{$\Lambda^\dagger$}
and fix a total order $>$ on $\Lambda^\dagger\simeq\mathbb N^{n}$.
We fix a $\Psi$-weighted lex order ``$>$'' on $\mathbb N^N$ and define a total
order ``$\triangleright$'' on $\Lambda^\dagger \times \mathbb N^N$ by: $(\lambda,\mathbf m)\triangleright(\mu,\mathbf m')$
if  $\lambda > \mu$, and if $\lambda = \mu$, then we set  $(\lambda,\mathbf m)\triangleright(\lambda,\mathbf m')$
if $\mathbf m>\mathbf m'$.
Given  $p(\mathbf x,\mathbf e^{\lambda})\in \mathbb C[{\cal Z}_S]$,
we define a $\Lambda\times\mathbb Z^N$-valued valuation:
\begin{equation}\label{valuationdef2}
\nu(p({\bf x}, e^{\lambda}))=\min\{(\lambda,{\mathbf p}) \in  \Lambda^\dagger\times\mathbb N^N\mid a_{\lambda,\mathbf p}\not=0\}
\end{equation}
for $p(\mathbf x, e^{\lambda}) = \sum_{(\lambda, {\bf p})\in  \Lambda^\dagger\times\mathbb N^N} a_{\lambda,\bf p} {\bf x}^{\bf p} e^\lambda$.
For a rational function $h = \frac{p}{p'}$ we define $\nu(h)=\nu(p)-\nu(p')$.
The valuation $\nu$ is called the {\it lowest term valuation} with respect to the
para\-meters $x_{{1}},\ldots, x_{N}, e^{\omega_1},\ldots,e^{\omega_{n'}},e^{\eta_1},\ldots,e^{\eta_{n-n'}}$ and the monomial order ``$\,\ge$''.

By the argument in the end of section \ref{birat}, it makes sense to view $\nu$ also as a valuation on $\mathbb C(G/\hskip -3.5pt/U)$.
We associate to $\mathbb C[G/\hskip -3.5pt/U]$ the {\it valuation monoid} ${\cal V}(G/\hskip -3.5pt/U)$\index{${\cal V}(G/\hskip -3.5pt/U)$, valuation monoid}:
\begin{equation}\label{valuationsemigroup}
{\cal V}(G/\hskip -3.5pt/U)={\cal V}(G/\hskip -3.5pt/U,\nu,>)
=\{\nu({p})\mid p\in \mathbb C[G/\hskip -3.5pt/U]-\{0\}\}
\subseteq \Lambda\times\mathbb Z^N.
\end{equation}
\begin{remark}\label{nueinsbewertung} It is easy to check that
the monoid ${\cal V}(G/\hskip -3.5pt/U)$ is independent of the choice of the total order on $\Lambda^\dagger$.
\end{remark}
\subsection{Filtrations and essential elements}\label{Filtessent} Let $\mathbb N^N$ be endowed with a $\Psi$-weighted lex order ``$>$''.
\begin{definition} Let $S=(\beta_1,\ldots, \beta_N)$ be a birational sequence for $U^-$.
An increasing sequence of subspaces $U(\mathfrak n^-)_{\bf \le m}$
$\mathbf m\in \mathbb N^N$, defined by
$$
U(\mathfrak n^-)_{\bf \le m}=\langle\mathbf f^{(\mathbf k)}=f_{\beta_1}^{(k_1)}\cdots f_{\beta_N}^{(k_N)}\mid \mathbf k\le \mathbf m \rangle
$$
is called a $(\mathbb N^N,>,S)$-{\it filtration} of $U(\mathfrak n^-)$.
The associated graded vector space is denoted by $U^{gr}(\mathfrak n^-)=
\bigoplus_{{\bf m}\in \mathbb N^N} U^{gr}(\mathfrak n^-)_{\bf m}$, where
$$
U^{gr}(\mathfrak n^-)_{\bf m}:=U(\mathfrak n^-)_{\bf \le m}/U(\mathfrak n^-)_{\bf < m}\hbox{\ and\ }
U(\mathfrak n^-)_{\bf < m}=\langle \mathbf f^{(\mathbf k)}\mid \mathbf k< \mathbf m \rangle.
$$
\end{definition}
\begin{definition}
A multi-exponent $\mathbf m\in \mathbb N^N$ is called {\it essential} for the $(\mathbb N^N,>,S)$-filtration if
$U^{gr}(\mathfrak n^-)_{\bf m}\not=0$.
Denote by $es(\mathfrak n^-)\subseteq  \mathbb N^N$\index{$es(\mathfrak n^-)$, essential monoid} the set of all essential multi-exponents.
\end{definition}

Given a $(\mathbb N^N,>,S)$-filtration on $U(\mathfrak n^-)$,
we get an induced filtration on $V(\lambda)$ as follows:
$$
V(\lambda)_{\bf \le m}:=U(\mathfrak n^-)_{\bf \le m}  v_\lambda, \quad V(\lambda)_{\bf < m}:=U(\mathfrak n^-)_{\bf <m}  v_\lambda\quad \forall \mathbf{m}\in \mathbb N^N.
$$
We set
$$
V^{gr}(\lambda)=\bigoplus_{{\bf m}\in \mathbb N^N} V(\lambda)_{\bf m},\hbox{\  where\ }V(\lambda)_{\bf m}= V(\lambda)_{\bf \le m}/V(\lambda)_{\bf <m}.
$$
Obviously we have  $\dim V(\lambda)_{\bf m}\le 1$.
\begin{definition}
A pair  $(\lambda,{\bf m})\in \Lambda^+\times \mathbb N^N$ is called {\it essential} for the $(\mathbb N^N,>,S)$-filtration of $V(\lambda)$ if $V(\lambda)_{\bf m}\not=0$.
If $(\lambda,\bf m)$ is essential, then $\mathbf f^{(\mathbf m)}v_\lambda$ is called an {\it essential vector} for $V(\lambda)$ and
$\bf m$ is called an {\it essential multi-exponent} for $V(\lambda)$.
Let  $es(\lambda)\subseteq \mathbb N^N$\index{$es(\lambda)$, set of essential multi-exponents} be the set of all essential multi-exponents for $V(\lambda)$.
\end{definition}
The term {\it essential} makes sense only in connection with a fixed $(\mathbb N^N,>,S)$-filtration on $U(\mathfrak n^-)$.
We often omit the reference to a filtration and assume tacitly that we have fixed one. One can show that the following sets are naturally endowed with the structure of a monoid:
$$
es(\mathfrak n^-)\subseteq \mathbb N^N, \Gamma(\lambda)=\bigcup_{n\in\mathbb N} n\times es(n\lambda)\subset \mathbb N\times \mathbb N^N\ \hbox{and\ }\
\Gamma=\bigcup_{\lambda\in\Lambda^+} \{\lambda\}\times es(\lambda) \subset \Lambda^+\times \mathbb N^N.
$$
\begin{definition}
We call $es(\mathfrak n^-)$ the {\it essential monoid} of $\mathfrak n^-$,
$\Gamma(\lambda)$ \index{$\Gamma(\lambda)$, essential monoid associated to $\lambda$} is called the {\it essential monoid} associated to $\lambda$, and $\Gamma$
is called the {\it global essential monoid}. The real cone ${\cal C}_{(>,S)}=\overline{\mathbb R_{\ge 0} \Gamma} \subset \Lambda_{\mathbb R}\times \mathbb R^N$\index{${\cal C}_{(>,S)}$, essential cone}
is called the {\it essential cone} associated to $(\mathbb N^N,>,S)$.
\end{definition}
\begin{remark} If $\Gamma$ is finitely generated, then $\mathbb R_{\ge 0} \Gamma$ is obviously already a closed subset.
\end{remark}

Let $p:\Lambda_{\mathbb R}\times \mathbb R^N\rightarrow \Lambda_{\mathbb R}$ be the projection onto the first component.
\begin{definition}
For $\lambda\in \Lambda^+_{\mathbb R}$, the {\it essential polytope} $P(\lambda)$\index{$P(\lambda)$, essential polytope} associated to $\lambda$ is defined as
$P(\lambda)=p^{-1}(\lambda)\cap {\cal C}_{(>,S)}$.
\end{definition}

\begin{remark}\label{changepsi}
Let $S=(\beta_1,\ldots,\beta_N)$ be a birational sequence and fix a $\Psi$-weighted lex order ``$>$''
on $\mathbb N^N$. If one replaces $\Psi$ by $k\Psi$ for some $k>0$, then the
associated global essential monoid does not change because the filtration
does not change. So one should think of the possible choices for $\Psi$
as rational points on the intersection $\cal W$ of a sphere with the positive quadrant
in $\mathbb R^N$.

Further, let $\tilde \Psi$ be an integral weight function such that  $\tilde \Psi(\mathbf m)$
depends only on the weight $\sum_{i=1}^N m_i\beta_i$. An example for such a function
is given in Example 9. If one replaces in the situation above
$\Psi$ by $\Psi +\tilde\Psi$, then the filtration may change but the associated
global essential monoid does not change, this can be proved using an easy weight argument.

In fact, we conjecture that $\cal W$ admits a finite triangulation such that the global essential monoid
stays constant for any choice of $\Psi$ in the interior of a simplex of the triangulation.
\end{remark}

\subsection{The dual essential basis}\label{gammagleichnu}
By construction,  $\mathbb B_\lambda=\{\mathbf f^{(\mathbf p)}
v_\lambda \mid {\bf p} \in es(\lambda) \}$ is a basis of $V(\lambda)$, let $\mathbb B_\lambda^*\subset V(\lambda)^*$
be the dual basis. For $\mathbf f^{(\mathbf p)}v_\lambda\in \mathbb B_\lambda$ denote by $\xi_{\lambda,\bf p}$ the corresponding dual element in $\mathbb B_\lambda^*$.
Keeping in mind the isomorphisms in \eqref{coordinateGnachU}, the set
$\mathbb B_\Gamma:=\{\xi_{\lambda,\bf p}\mid (\lambda,{\bf p})\in \Gamma\}$\index{$\mathbb B_\Gamma$, dual essential basis}
is a vector space basis for the algebra $\mathbb C[G/\hskip -3.5pt/U]$, we call it the {\it dual essential basis}.
Consider the structure constants $c_{\lambda,{\bf p};\mu,\mathbf q}^{\lambda+\mu,\mathbf r}$, defined for $\mathbf p\in {\rm es}(\lambda)$ and
$\mathbf q\in {\rm es}(\mu)$ by
\[
\xi_{\lambda,\mathbf p}\xi_{\mu,\mathbf q}=\sum_{\mathbf r\in {\rm es}(\lambda+\mu)} c_{\lambda,\mathbf p;\mu,\mathbf q}^{\lambda+\mu,\mathbf r}
\xi_{\lambda+\mu,\mathbf r}.
\]
These basis elements satisfy the following multiplicative property:
$$\xi_{\lambda,\mathbf p}\xi_{\mu,\mathbf q}= \xi_{\lambda+\mu,\mathbf p+\mathbf q} + \sum_{\mathbf r\in {\rm es}(\lambda+\mu),\mathbf r >\mathbf p+\mathbf q}
c_{\lambda,\mathbf p;\mu,\mathbf q}^{\lambda+\mu,\mathbf r}\xi_{\lambda+\mu,\mathbf r}.$$
As an immediate consequence one can show \cite{FaFL}:
\begin{proposition}\label{semigleichsemi1}\it
The valuation {monoid} ${\cal V}(G/\hskip -3.5pt/U)$ and the global essential monoid $\Gamma$ coincide.
\end{proposition}

\subsection{ASM-sequence and a toric variety}
Let $\mathbb T$ be the torus $(\mathbb C^*)^N$.
The algebra $\mathbb C[\Gamma]$ associated to the {monoid} $\Gamma$
can be naturally endowed with the structure of a $T\times \mathbb T$-algebra by
\begin{equation}\label{TACTION}
(t,t_1,\dots,t_N)\cdot (\lambda,\mathbf p):=\lambda(t)\big(\prod_{i=1}^N t_i^{p_i}\big)(\lambda,\mathbf p).
\end{equation}
\begin{definition}
Let $S=(\beta_1,\ldots,\beta_N)$ be a birational sequence and fix a $\Psi$-weighted lex order ``$>$'' on $\mathbb N^N$.
We call $(S,>)$\index{$(S,>)$, a sequence with an ASM} a {\it sequence with an affine, saturated monoid} (short ASM) \index{$ASM$, affine saturated monoid}
if $\Gamma$ is finitely generated and saturated.
\end{definition}
The two assumptions imply that (see \cite{CLS} for more details on toric varieties) $\hbox{\textrm{Spec}\,} (\mathbb C[\Gamma])$ is naturally
endowed with the structure of a normal toric variety for $T\times \mathbb T$ of dimension $\dim T+\dim \mathfrak n^-$.

\begin{remark}
If $\lambda\in\Lambda^+$, then $(S,>)$ being a sequence with an ASM implies that the lattice points of the essential polytope $P(\lambda)$ are exactly the pairs $(\lambda,\mathbf m)$ such that $\mathbf m$ is an essential multi-exponent for $\lambda$.
\end{remark}

\subsection{A filtration of \texorpdfstring{$\mathbb C[G\hbox{/\hskip -3.5pt/U]}$}{the coordinate ring of the flag variety}} \label{algfilt}
In the following set $R=\mathbb C[G{/\hskip -3.5pt/}U]$.
Let ``$>$" be the fixed $\Psi$-weighted lex order on $\mathbb N^N$. We
define a new partial order ``$>_{alg}$" on $\Lambda^+\times\mathbb N^N$ as follows:
\begin{definition}\label{algorder}
$(\lambda,\mathbf p)>_{alg}(\mu,\mathbf q)\hbox{\ if\ }\lambda\succ_{wt} \mu\hbox{\ or\ }\lambda=\mu \hbox{\ and\ }\mathbf p<\mathbf q$.
\end{definition}
Note that we turn around ``$>$" when going to the coordinate ring.
We define a filtration as follows: for $\lambda\in\Lambda^+$ and $\mathbf p\in es(\lambda)$ set
\begin{equation}\label{deg1}
R_{\le_{alg}(\lambda,\mathbf p)}
=\bigoplus_{\substack{(\mu,\mathbf q)\in\Gamma\\ (\mu,\mathbf q)\le_{alg}(\lambda,\mathbf p)}} \mathbb C\xi_{\mu,\mathbf q}
\hbox{\ and\ }
R_{<_{alg}(\lambda,\mathbf p)}=\bigoplus_{\substack{(\mu,\mathbf q)\in\Gamma\\(\mu,\mathbf q)<_{alg}(\lambda,\mathbf p)}} \mathbb C \xi_{\mu,\mathbf q}.
\end{equation}
The associated graded vector space $gr R$ is again a ring and isomorphic to $\mathbb C[\Gamma]$.
\par
Now assume that $(S,>)$ a sequence with an ASM. As in \cite{AB}, one can change the multi-filtration  into an
$\mathbb N$-filtration. Translated into geometry one can reformulate this into: the affine $G$-algebra $R$ corresponds to
the affine $G$-scheme $Y= \hbox{\textrm{Spec}}(R)$, and $gr R$ corresponds to an affine
$T\times \mathbb T$-scheme denoted by $Y_0$, and there exists a family of affine $T$-schemes
$\rho : \mathcal Y\rightarrow \mathbb A^1$ such that

\begin{proposition}\label{Y0}\it
\begin{enumerate}
\item  $\rho$ is flat;
\item $\rho$ is trivial with fiber $Y$ over the complement of $0$ in $\mathbb A^1$;
\item and
the fiber of $\rho$ at $0$ is isomorphic to $Y_0$.
\end{enumerate}
\end{proposition}
\subsection{Degeneration of polarized projective \texorpdfstring{$G$}{G}-varieties}
We follow again the approach of Alexeev and Brion \cite{AB}.
Let $(Y,\cal L)$ be a  polarized projective $G$-variety, \emph{i.e.}, $Y$ is a normal projective $G$-variety together with an
ample $G$-linearized  invertible sheaf $\cal L$.

The sheaf ${\cal L}^n := {\cal L}^{\otimes n}$ is also $G$-linearized, so the space
of global sections $H^0(Y,{\cal L}^n)$ is a finite-dimensional rational $G$-module.
Consider the associated graded algebra
$$
R(Y,{\cal L}):= \bigoplus_{n\in\mathbb N} H^0(Y,{\cal L}^n),
$$
this is a finitely generated, integrally closed domain. We have $Y = \hbox{Proj\,} R(Y,{\cal L})$ and
${\cal L}^n={\cal O}_Y(n)$.

One can endow the algebra $R(Y,{\cal L})$ with a filtration as in the section before.
The associated algebra $gr R(Y,{\cal L})$ is still finitely generated, it is an
integrally closed domain, with an action of $\mathbb C^*\times T \times \mathbb T$
such that the $\mathbb C^*$-action defines a positive grading.
The projective variety $Y_0 :=\hbox{Proj\,} gr R(X,\cal L)$ is again
a projective $T\times\mathbb  T$-variety equipped with $T\times\mathbb  T$-linearized sheaves
${\cal L}_0^{(n)}={\cal O}_{Y_0}(n D)$ for all integers $n$, where $D$ refers to a $\mathbb Q$-Weil divisor.
Further, $D$ is $\mathbb Q$-Cartier and ample, \emph{i.e.}, the sheaf ${\cal L}_0^{(m)}$
 is invertible and ample for any sufficiently divisible integer $m > 0$.
In particular, every sheaf ${\cal L}_0^{(n)}$ is divisorial, \emph{i.e.}, it is the sheaf of sections of an integral Weil divisor.
We just quote \cite{AB}, the proof is the same:
\begin{theorem}\it
Let $(Y,\cal L)$ be a  polarized $G$-variety and let
$(S,>)$ be a sequence with an ASM. Consider the
induced $(\mathbb N^N,>,S)$-{\it filtration} on $R(Y,{\cal L})$.
Then there exists a family of $T$-varieties $\pi : {\cal Y}\rightarrow \mathbb A^1$, where ${\cal Y}$ is a
normal variety, together with divisorial sheaves ${\cal O}_{\cal Y}(n)$ ($n\in\mathbb  Z$), such that
\begin{itemize}
\item[{\it i)}] $\pi$ is projective and flat.
\item[{\it ii)}] $\pi$ is trivial with fiber $Y$ over the complement of $0$ in $\mathbb A^1$, and ${\cal O}_{\cal Y}(n)\vert_Y \simeq {\cal L}^n$
for all $n$.
\item[{\it iii)}]  The fiber of $\pi$ at $0$ is isomorphic to $Y_0$, and ${\cal O}_{\cal Y}(n)\vert_{Y_0} \simeq {\cal L}_0^{(n)}$ for all $n$.
\end{itemize}
In addition, if $Y$ is spherical, then $Y_0$ is a toric variety for the torus $T\times\mathbb T$.
\end{theorem}
\subsection{Moment polytopes}
The results and definitions in this section on moment polytopes can be found in \cite{AB}, the only
difference being that we consider not only the string cone case.
Let $(Y,\cal L)$ be a  polarized spherical $G$-variety and let
$(S,>)$ be a sequence with an ASM,
so $\Gamma$ is assumed to be finitely generated and saturated.
We recall the definition of a moment polytope of the polarized $G$-variety $(Y,\cal L)$.
Note that for a dominant weight $\lambda$, the isotypical component $H^0(Y ,{\cal L}^{n})_{(\lambda)} \not= 0$
if and only if the space of $U$ invariants of weight $\lambda$ is not trivial: $H^0(Y,{\cal L}^n)^U_\lambda \not=0$.
Further, the algebra $R(Y,{\cal L})^U$ is finitely generated and $\Lambda^+\times\mathbb N$-graded;
let $(f_i)_{i=1,\ldots,r}$ be homogeneous
generators and $(\lambda_i, n_i)$ their weights and degrees.

\begin{definition}
The convex hull  of the points $\frac{\lambda_i}{n_i}$, $i=1,\ldots,r$, in $\Lambda_{\mathbb R}$
is called the {\it moment polytope $P(Y,\cal L)$}\index{$P(Y,\cal L)$, moment polytope} of the polarized $G$-variety $(Y,\cal L)$.
\end{definition}
Another way to view the moment polytope is the following: the points $\frac{\lambda}{n}\in\Lambda_{\mathbb R}$
such that $\lambda\in\Lambda^+$, $n\in\mathbb N_{>0}$, and the isotypical component $H^0(Y,{\cal L}^n)_{(\lambda)}$
is nonzero, are exactly the rational points of the rational convex polytope $P(Y,{\cal L})\subset \Lambda_{\mathbb R}$.
Further,  $P(Y,{\cal L}^m) = mP(Y,{\cal L})$ for any positive integer $m$.

By positive homogeneity of the moment polytope, this definition extends to
$\mathbb Q$-polarized varieties, in particular, to any limit $(Y_0,{\cal L}_0)$ of $(Y,{\cal L})$.
We denote by $P(\Gamma,Y,{\cal L})$ the moment polytope of that limit. It is a rational convex polytope in
$\Lambda_{\mathbb R}\times \mathbb R^N$, related to the moment polytope of $(Y,{\cal L})$ by the following theorem:

\begin{theorem}\label{polytopethm}\it
The projection $p : \Lambda_{\mathbb R}\times \mathbb R^N \rightarrow \Lambda_{\mathbb R}$ onto the first component
restricts to a surjective map
$$
p : P(\Gamma,Y,{\cal L}) \rightarrow P(Y,{\cal L}),
$$
with fiber over any $\lambda\in\Lambda^+_{\mathbb R}$
being the essential polytope $P(\lambda)$. In particular, for $\lambda\in\Lambda^+$, let
$P_\lambda$ be the stabilizer  of the line $[v_\lambda]\in \mathbb P(V(\lambda))$.
The limit of the flag variety $G/P_\lambda$ is a toric variety under $\mathbb T$, and its
moment polytope is the essential polytope $P(\lambda)$.
\end{theorem}
\subsection{Essential polytope and Newton-Okounkov body}\label{NObody}
If we restrict the valuation $\nu$ to the subring $\mathbb C[x_1,\ldots,x_N]$
by omitting the dominant weights in the definition in \eqref{valuationdef2},
then we get a $\mathbb Z^N$-valued valuation $\nu_1$, which, by Remark~\ref{nueinsbewertung},
more or less completely determines the valuation monoid.

The valuation $\nu_1$ depends on the choice of $S$  and ``$>$'', and the birational map
$\pi:Z_S\rightarrow U^-$ (see \eqref{birational1}) provides a birational map $\pi:Z_S\rightarrow G/B$ by identifying $U^-$ with an affine neighborhood
of the identity. Hence we can view $\nu_1$ naturally as a $\mathbb Z^N$-valued valuation on $\mathbb C(G/B)$.

Let $\lambda$ be a regular dominant weight and let $\mathcal L_\lambda$
be the corresponding very ample line bundle on $G/B$. Let $R_\lambda$ be the ring $\bigoplus_{n\ge 0} H^0(G/B,\mathcal L_{n\lambda})$.
Recall that $H^0(G/B,\mathcal L_{n\lambda})\simeq V(n\lambda)^*$ as a $G$-representation, we fix for all $n\in\mathbb N$ the
dual vector $\xi_{n\lambda,\mathbb O}$ to the fixed highest weight vector $v_{n\lambda}\in V(n\lambda)$. The  Newton-Okounkov body
associated to $G/B$ depends on the choice of the valuation $\nu_1$, the ample line bundle $\mathcal L_\lambda$ and the
choice of a non-zero element in $H^0(G/B,\mathcal L_{\lambda})$, in our case the vector $\xi_{n\lambda,\mathbb O}$.
The {\it Newton-Okounkov body} $\Delta_{\nu_1}(\lambda)$\index{$\Delta_{\nu_1}(\lambda)$, Newton-Okounkov body} is defined as follows.
(For more details on Newton-Okounkov bodies see for example \cite{KK,K1}). One associates to the graded ring $R_\lambda$ the monoid
$$
\bigcup_{n>0}\{(n,\nu_1(\frac{s}{\xi_{n\lambda,\mathbb O}})) \mid s\in H^0(G/B,\mathcal L_{n\lambda})\}.
$$
In view of Proposition~\ref{semigleichsemi1}, this is the essential monoid $\Gamma(\lambda)$ associated to $\lambda$.
\begin{definition}\label{Def:NOBody}
The {\it Newton-Okounkov body} $\Delta_{\nu_1}(\lambda)$ is the convex closure
$$
\textrm{conv}(\overline{\bigcup_{n\in \mathbb N}\{\frac{\mathbf m}{n}\mid   \mathbf m\in es(n\lambda)\}}).
$$
\end{definition}
It follows immediately that the Newton-Okounkov body $\Delta_{\nu_1}(\lambda)$ coincides with the essential polytope $P(\lambda)$.

\section{Examples of saturated, finitely generated monoids $\Gamma$}\label{examplesectionASM}
\subsection{Reduced decompositions and string polytopes}\label{Sec:String}
Fix a reduced decomposition $\underline{w}_0=s_{i_1}\cdots s_{i_N}$ of the longest word in the Weyl group of $\mathfrak g$
and let ${C}_{\underline{w}_0}\subset \mathbb R^N$\index{${C}_{\underline{w}_0}$, string cone} be the associated string cone defined in \cite{BZ,L3} (see Example \ref{Exo:String} below for an explicit description for $G=SL_n$)
and set $S=({\alpha_{i_1}},\ldots, {\alpha_{i_N}})$. By Example~\ref{monomialszwo}, we know that
$S$ is a birational sequence. Let $\Psi:\mathbb N^N\rightarrow \mathbb N$ be the height weighted function
as in Example~\ref{rootorder}. We fix on $\mathbb N^N$ the associated $\Psi$-weighted opposite lexicographic order.
For an element $\mathbf m\in {C}_{\underline{w}_0}$ denote by $G(\mathbf m)$ the corresponding element
of the global crystal basis of $U(\mathfrak n^-)$ \cite{Ka}, specialized at $q=1$.
The cone ${\cal C}_{\underline{w}_0}\subset\Lambda_{\mathbb R}\times\mathbb R^N$\index{${\cal C}_{\underline{w}_0}$, global string monoid} (see \cite{BZ,L3}, compare also \cite{AB})
is defined to be the intersection of $\Lambda_{\mathbb R}\times {C}_{\underline{w}_0}$ with $N$ half-spaces:
\begin{equation}\label{weightequation}
{\cal C}_{\underline{w}_0}=\left\{(\lambda,\mathbf m)\in \Lambda_{\mathbb R}\times {C}_{\underline{w}_0}\mid
\begin{array}{l}
m_k \le \langle\lambda,\alpha_{i_k}^\vee\rangle -\sum_{\ell=k+1}^N\langle \alpha_{i_\ell},\alpha_{i_k}^\vee \rangle m_\ell,\\ 
k=1,\ldots,N
\end{array}\right\}.
\end{equation}
\begin{theorem}[\cite{FaFL}]\label{stringCone}
We have $es(\mathfrak n^-)={C}_{\underline{w}_0}\cap \mathbb Z^N$.
\end{theorem}
As a corollary, $(S,>)$ is a sequence with an ASM. Moreover, the proof of the theorem shows more precisely that the $(\mathbb N^N,>,S)$-filtration of $U(\mathfrak n^-)$ is compatible with Kashiwara's global crystal basis. In particular, given a dominant weight $\lambda$, $\mathbf m\in   \mathbb N^N$ is an essential multi-exponent for $V(\lambda)$ if and only if $G(\mathbf m).v_\lambda\not=0$.
\par
For $\lambda\in\Lambda^+$, the essential polytope $P(\lambda)$ is called the \emph{string polytope} and will be denoted by $Q_{\underline{w}_0}(\lambda)$\index{$Q_{\underline{w}_0}(\lambda)$, string polytope}.

\begin{exam}\label{Exo:String}\rm
We give a description of the string cone for $G=SL_n$. We follow here the description of \cite{BZ}, while a recursive one can be found in \cite{L3}. For any fundamental weight $\omega_i$, let
$w^i$ be the minimal representative of the coset $W_i s_i w_0$ where $W_i$ is the maximal parabolic subgroup generated by $s_j$ for $j \neq i$.
\par
Fix a reduced decomposition ${\underline{w}_0} = s_{i_1} \cdots s_{i_N}$ and $1 \leq i \leq n-1$. Let $\mathbf{s} = s_{i_{k(1)}} \cdots s_{i_{k(p)}}$ be a subword of ${\underline{w}_0} $ which is a reduced word of $w^i$. We set $T(\mathbf{s})$ to be the half-space defined by the inequality:
\[
\sum_{j=0}^p \sum_{k(j) < k < k(j+1)} \langle \omega_i,s_{i_{k(1)}} \cdots s_{i_{k(j)}}(\alpha_{i_k}^\vee)\rangle x_k \geq 0,
\]
Then the string cone $C_{{\underline{w}_0}} \subset \mathbb{R}^N_{\geq 0}$ is the intersection of all $T(\mathbf{s})$, where $\mathbf{s}$ is a subword of ${\underline{w}_0}$ and a reduced word for some $w^i$.
\end{exam}
\begin{example}\rm
We recall here the example 3.15 in \cite{BZ} for $G = SL_4$. The string cone $C_{s_2s_1s_3s_2s_1s_3} \subset \mathbb{R}^6_{\geq 0} $ is defined by
\[
x_4 - x_5 - x_6 \geq 0 \, , \, x_3 - x_5 \geq 0 \, ,\,x_2 - x_6 \geq 0 \, ,\, x_2 + x_3 - x_4 \geq 0.
\]
To obtain the global essential monoid, we have to determine the string polytopes and cut out by the weight-restrictions (\ref{weightequation}). 
For $\lambda = m_1\omega_1 + m_2 \omega_2 + m_3 \omega_3 \in \Lambda_{\mathbb{R}}^+$, the string polytope $Q_{{\underline{w}_0}}(\lambda)$ 
is the intersection of  $C_{s_2s_1s_3s_2s_1s_3} $ with the half-spaces defined by the inequalities:
$$
\begin{array}{c}
x_6 \leq m_3 \, , \, x_5 \leq m_1  \, , \, x_4 \leq m_2 + m_1 + m_3  \, , \\ \,  x_3 \leq m_2 - m_3  \, , \,  x_2 \leq m_2 - m_1  \, , \,  x_1 \leq 2m_1 + 2m_3 - m_2.
\end{array}
$$
The cone $\mathcal{C}_{s_2s_1s_3s_2s_1s_3}$ is then
\[
\mathcal{C}_{s_2s_1s_3s_2s_1s_3} = \bigcup_{\lambda \in \Lambda_{\mathbb{R}}^+} (\lambda, Q_{{\underline{w}_0}}(\lambda) ).
\]
\end{example}

\begin{remark}
For a fixed $\lambda\in\Lambda^+$, let $Q(\lambda)$ be the set of string polytopes $Q_{\underline{w}_0}(\lambda)$ when $\underline{w}_0$ runs through all reduced decompositions of $w_0$. An equivalent relation given by the unimodularly equivalence can be defined on $Q(\lambda)$. It is interesting to find representatives for the equivalent classes. For example, it is clear that if $\underline{w}_0'$ is obtained from ${\underline{w}_0}$ by an exchange $s_i s_j \mapsto s_j s_i$ where $s_i, s_j$ are orthogonal reflections, $Q_{{\underline{w}_0}}(\lambda)$ is unimodularly equivalent to $Q_{{\underline{w}_0'}}(\lambda)$. But this condition is not enough to determine the equivalence classes.
\end{remark}

\subsection{The Gelfand-Tsetlin polytopes}
We recall the most famous example of toric degenerations of flag varieties, namely the toric degeneration arising from Gelfand-Tsetlin polytopes \cite{GT}, and explain how these polytopes are related to the string polytopes. We fix $N=\binom{n}{2}$ and $M=\binom{n+1}{2}$.
\begin{definition}\label{gt-polytope}
Let $\lambda = \sum_{i =1}^{n-1} m_i \omega_i \in \Lambda^+$ be a dominant integral weight of $G=SL_n$ and $m_n=0$. The Gelfand-Tsetlin polytope associated to $\lambda$ is defined as
\[
GT(\lambda)\index{$GT(\lambda)$, Gelfand-Tsetlin polytope} :=  \left\{ (x_{i,j})_{0 \leq i \leq n-1, 1 \leq j \leq n - i} \in \mathbb{R}^M \mid \begin{array}{c} x_{0,j} = \sum_{i = j}^n m_i \\   x_{i,j} \geq x_{i+1, j}  \, , \, x_{i,j} \geq  x_{i-1, j+1} \end{array} \right\}.
\]
The set of lattice points in $GT(\lambda)$ is denoted by:
\[
GT_{\mathbb{Z}}(\lambda)\index{$GT_{\mathbb{Z}}(\lambda)$, lattice points in Gelfand-Tsetlin polytope} = GT(\lambda) \cap \mathbb{Z}^M.
\]
\end{definition}
\begin{example}
Let $G = SL_4$ and $\lambda = (\lambda_1-\lambda_2) \omega_1 + (\lambda_2-\lambda_3) \omega_2 + \lambda_3 \omega_3\in \Lambda^+$. Then
\[
GT(\lambda) :=  \left\{ 
{\tiny
(x_{i,j})_{0 \leq i \leq 3, 1 \leq j \leq 4 - i} \in \mathbb{R}^{10} \mid
\begin{array}{c} x_{0,1} = \lambda_1, x_{0,2} =\lambda_2, x_{0,3} = \lambda_3, x_{0,4} =  0\\  x_{0,1} \geq x_{1,1} \geq x_{0,2} \geq x_{1,2} \geq x_{0,3} \geq x_{1,3} \geq x_{0,4} \\
x_{1,1} \geq x_{2,1} \geq x_{1,2} \geq x_{2,2} \geq x_{1,3} \\
x_{2,1} \geq x_{3,1} \geq x_{2,2}
\end{array} 
}
\right\}.
\]
\end{example}
\begin{proposition}[\cite{GT}]
\begin{enumerate}
\item For any $\lambda, \mu \in \Lambda^+$, we have:
\[
GT_{\mathbb{Z}}(\lambda)  + GT_{\mathbb{Z}}(\mu)  = GT_{\mathbb{Z}}(\lambda + \mu).
\]
\item The monoid $\Gamma_{GT}\index{$\Gamma_{GT}$, Gelfand-Tsetlin monoid} := \bigcup_{\lambda \in \Lambda^+} (\lambda, GT_{\mathbb{Z}}(\lambda) )$ is finitely generated and saturated.
\end{enumerate}
\end{proposition}

\begin{remark}
\begin{enumerate}
\item The toric degeneration constructed in \cite{GL} is in fact isomorphic to the one associated with $\Gamma_{GT}$ \cite{KM}.
\item Let ${\underline{w}_0} = s_1\, (s_2s_1) \, (s_3s_2s_1) \, \ldots \, (s_{n-1}s_{n-2} \ldots s_1)$ and $\lambda \in \Lambda^+$. The string polytope $Q_{{\underline{w}_0}}(\lambda)$ is unimodularly equivalent to $GT(\lambda)$ \cite{L3}.
\end{enumerate}
\end{remark}

It has been shown in \cite{KM} that any Schubert variety degenerates to a union of toric varieties associated to certain faces of the Gelfand-Tsetlin polytope and they have given the following formula in terms of dual Kogan faces of the polytope. The following description is due to \cite{KST}:
\par
A \textit{dual Kogan facet} is a facet  of $GT(\lambda)$
 defined by an equality of the form $x_{i,j} = x_{i+1, j-1}$. We associate to such a facet the simple reflection $s_{n+1 - j}$.
\par
Suppose now $F_\lambda$ is a face of $GT(\lambda)$ defined by the intersection of dual Kogan facets. Then we compose the simple reflections corresponding to the equalities by following $i =n, \ldots, 1$ (bottom to top) and within this $j = n+1 - j , \ldots , 1$ (right to left). The resulting Weyl group element will be denoted by $w(F_\lambda)$.
\begin{theorem}[\cite{KM}]
The Schubert variety $X_w$ degenerates to the reduced union of toric varieties associated to the face $F_\lambda$ satisfying $w(F_\lambda) = ww_0$.
\end{theorem}
\begin{remark}\label{rem-kempf}
If $w$ is a Kempf element (see \cite{La} for the definition), then the toric variety is irreducible \cite{KST}. In \cite{DY2}, a 
construction of toric degenerations of Schubert varieties has been studied also for more general Weyl group elements.
\end{remark}

\begin{remark}
Let $w \in W$ and fix a reduced decomposition $\underline{w} = s_{i_1} \cdots s_{i_s}$. We can extend $\underline{w}$ to a reduced decomposition ${\underline{w}_0} = s_{i_1} \cdots s_{i_s} s_{i_{s+1}} \cdots s_{i_N}$. We denote by $F_{\underline{w}} \subset Q_{{\underline{w}_0}}(\lambda)$ the face of the string polytope defined by setting the coordinates $x_{i_{s+1}}, \ldots, x_{i_N}$ to $0$. Then the irreducible toric variety associated to the face $F_{\underline{w}}$ gives a toric degeneration of the Schubert variety $X_w$.
\par
We see that in the case of string polytopes, we can always \text{choose} for a particular Schubert variety $X_w$ a toric degeneration of the whole flag variety being compatible with $X_w$. So there is always an irreducible toric degeneration. Compare this with the case of the Gelfand-Tsetlin degeneration discussed, where we can find these toric degenerations of the Schubert varieties within a fixed toric degeneration of $G/P_\lambda$ but at the price of this variety being not irreducible.
\end{remark}
\begin{remark}
For a general discussion about degenerations of Schubert varieties, Richardson varieties and Kazhdan-Lusztig varieties see \cite{Kn}.
\end{remark}

\subsection{Birational sequences arising from convex orders and Lusztig polytopes}\label{Lusztig}
Fix a reduced decomposition $\underline{w}_0=s_{i_1}\cdots s_{i_N}$ of the longest word $w_0$ in the Weyl group and let $S=(\beta_1,\ldots,\beta_N)$ be an enumeration of the positive roots associated to the decomposition, \emph{i.e.}, $\beta_k=s_{i_1}\cdots s_{i_{k-1}}(\alpha_{i_k})$ for $k=1,\ldots,N$. By Example~\ref{monomialseins} we know that
$S$ is a birational sequence. Let $\Psi:\mathbb N^N\rightarrow \mathbb N$ be the height weighted function
as in Example~\ref{rootorder}. We fix on $\mathbb N^N$ the associated  $\Psi$-weighted opposite right
 lexicographic order ``$>$''.
For an element $\mathbf m\in \mathbb N^N$ denote by $B(\mathbf m)$ the corresponding element
of Lusztig's canonical basis of $U(\mathfrak n^-)$ (using Lusztig's parametrization with respect to the fixed decomposition $\underline{w}_0$),
specialized at $q=1$.

The decomposition determines also an enumeration of the positive roots (see the choice of $S$ above) and hence a PBW-basis.
The connection between this PBW-basis and the canonical basis was described by Lusztig in the $\tt A$, $\tt D$ and $\tt E$-case
\cite{Lu1,Lu2}, and by Caldero  \cite{Ca1} in the semisimple case. In terms of the $(\mathbb N^N,>,S)$-filtration
this can be reformulated as follows \cite{FaFL}:
\begin{theorem}\it
\begin{enumerate}
\item The $(\mathbb N^N,>,S)$-filtration of $U(\mathfrak n^-)$ is compatible with Lusztig's canonical basis.
In particular, given a dominant weight $\lambda$, $\mathbf m\in   \mathbb N^N$ is an essential
multi-exponent for $V(\lambda)$ if and only if $B(\mathbf m).v_\lambda\not=0$.
\item $(S,>)$ is a sequence with an ASM.
\end{enumerate}
\end{theorem}
So the convex hull $\mathcal{L}_S(\lambda)$\index{$\mathcal{L}_S(\lambda)$, Lusztig polytope} of the points $\{\mathbf m\in\mathbb N^N\mid B(\mathbf m).v_\lambda\not=0\}$
is the essential polytope, and the positive real span of $\Gamma_\mathcal{L}=\bigcup_{\lambda\in\Lambda^+} (\lambda, \mathcal{L}_S(\lambda))$\index{$\Gamma_\mathcal{L}$, Lusztig monoid}
is the essential cone ${\cal C}_{(>,S)}$. Since the construction uses Lusztig's parametrization of the canonical basis,
we call ${\cal C}_{(>,S)}$ the Lusztig cone and refer to $\mathcal{L}_S(\lambda)$ as the Lusztig polytope associated to the reduced 
decomposition  $\underline{w}_0$.

\begin{example}\rm
Let $\underline{w}_0$ be a reduced decomposition of $w_0$ and $S=(\beta_1,\ldots,\beta_N)$ be the sequence of positive roots associated to the reduced decomposition. A unimodular equivalence between the string polytope $Q_{\underline{w}_0}(\lambda)$ and the Lusztig polytope $\mathcal{L}_S(\lambda)$ is given in \cite{MG}, Corollaire 3.5.
\end{example}



\subsection{Birational sequence arising from good sequences and PBW polytopes}\label{GoodOrdering}
Let $G=SL_n$. In \cite{FFL2}, the PBW filtration on finite-dimensional simple modules was studied. For any simple module a polytope $\mathcal{P}(\lambda)$ was introduced whose lattice points parametrize a monomial basis of the associated graded module. These polytopes had been suggested by E.~Vinberg.
\par
We briefly recall the definition here. A sequence ${\bf b}=(\delta_1,\dots,\delta_r)$ of positive roots is called a {\it Dyck path} if the first and the last roots are simple roots ($\delta_1=\alpha_{i,i}=\alpha_i$, $\delta_r=\alpha_{j,j}=\alpha_j$, $i\le j$), and if $\delta_m=\alpha_{p,q}$, then $\delta_{m+1}$ is either $\alpha_{p+1,q}$ or $\alpha_{p,q+1}$.
\par
For $\lambda \in \Lambda^+$, we define a polytope
\begin{equation}\label{FFLpoly}
\mathcal{P}(\lambda)\index{$\mathcal{P}(\lambda)$, polytope arising from PBW filtration} := \left\{  (x_\alpha) \in  \mathbb R_{\ge 0}^{\Phi^+} \mid {\small
\begin{array}{c} \forall i=1,\ldots,n-1,\ \forall \text{ Dyck paths } {\bf b}=(\delta_1,\dots,\delta_r)\\ 
\text{ starting in }\alpha_i, \text{ ending in }\alpha_j:\\ 
\sum_{\ell=1}^r x_{\delta_\ell}\le (\lambda, \alpha_i + \ldots + \alpha_j)
\end{array}}
\right\}.
\end{equation}
The set of its lattice points and the corresponding monoid are defined by:
$$
\mathcal{P}_{\mathbb{Z}}(\lambda)\index{$\mathcal{P}_{\mathbb{Z}}(\lambda)$, lattice points in $\mathcal{P}(\lambda)$} = \mathcal{P}(\lambda)  \cap \mathbb{Z}^{\Phi^+} \; \text{ and } \; \Gamma_\mathcal{P}\index{$\Gamma_\mathcal{P}$, the monoid arising from PBW filtration} = \bigcup_{\lambda \in \Lambda^+} (\lambda, \mathcal{P}_{\mathbb{Z}}(\lambda)).
$$

\begin{theorem}[\cite{FFL2}]
For $\lambda \in \Lambda^+$, $\mathcal{P}_{\mathbb{Z}}(\lambda)$ parametrizes a monomial basis of the finite dimensional simple module of highest weight $\lambda$. Moreover, $\Gamma_\mathcal{P}$ is finitely generated and saturated.
\end{theorem}
\begin{example}\label{Ex:FFLSL4}
Let $G = SL_4$, $\lambda = m_1 \omega_1 + m_2 \omega_2 + m_3 \omega_3\in\Lambda^+$ and $(x_1,\cdots,x_6)$ be the coordinates corresponding to positive roots $(\alpha_1, \alpha_{1} + \alpha_2, \alpha_2 , \alpha_1 + \alpha_2 +\alpha_3 , \alpha_2+\alpha_3, \alpha_3)$. The polytope $\mathcal{P}(\lambda)$ is given by:
$$
\mathcal{P}(\lambda) =  \left\{  (x_1,\ldots,x_6) \in  \mathbb R_{\ge 0}^6 \mid 
{\tiny
\begin{array}{c} x_1 \leq m_1 \, , \, x_3 \leq m_2 \, ,\, x_6 \leq m_3  \, ,  \\  x_1 + x_2 + x_3 \leq m_1 + m_2 \, , \, x_3  + x_5 + x_6 \leq m_2 + m_3  \, , \\ x_1 + x_2  + x_4 + x_5  + x_6 \leq m_1 + m_2  + m_3 , \\ x_1 + x_2 + x_4  + x_5  + x_6 \leq m_1  + m_2  + m_3.
\end{array}
}\right\}
$$
\end{example}

We fix an ordering on the set of positive roots $\beta_1,\dots,\beta_N$ and assume that
\[
\beta_i\succ_{wt} \beta_j \text{ implies } i<j.
\]
An ordering with this property (the larger roots come first) will be called a {\it good ordering}.
Once we fix such a good ordering, this induces an ordering on the basis vectors $f_\beta$.
\par
We consider a sequence $S=(\beta_1,\beta_2,\cdots,\beta_N)$ arising from a good ordering. 
By Example~\ref{monomialseins} we know that
$S$ is a birational sequence. Let $\Psi:\mathbb N^N\rightarrow \mathbb N$ be the homogeneous degree function
as in Example~\ref{homorder}. We fix on $\mathbb N^N$ the associated  $\Psi$-weighted right 
lexicographic order "$>$''. Then the monoid $\Gamma_\mathcal{P}$ is the global essential monoid corresponding to $(S, >)$.
\par
Note that this description is independent of the choice of the ordering of the roots, as long as the ordering is a ``{\it good\/}'' ordering.

\begin{remark}\label{remark-pbw}
\begin{enumerate}
\item It has been shown in \cite{Fo1} that $\mathcal{P}(\lambda)$ is not unimodularly equivalent to $GT(\lambda)$ as long as 
$\lambda$ is regular and the rank of the Lie algebra is greater or equal to $3$.
\item It has been shown in \cite{FFL2} that the monoid $\Gamma$ is in fact generated by the union of 
$(\omega_i, \mathcal{P}_{\mathbb{Z}}(\omega_i))$.
\item Moreover, one can show using \textit{polymake} \cite{GJ}, that for $G = SL_6$ there exists a 
dominant weight $\lambda \in \Lambda^+$ such that 
$\mathcal{P}(\lambda)$ is not unimodularly equivalent to any string polytope $Q_{{\underline{w}_0}}(\lambda)$.
\item For regular $\lambda \in \Lambda^+$, the number of facets of $\mathcal{P}(\lambda)$ is given by 
$$
\frac{n(n-1)}{2} + \sum_{i=1}^{n-1} {i} C_{n-1-i},
$$ 
where $C_j$ is the Catalan number \cite{Fo1}.
\item 
Let $G=SL_n$ and $w:\Phi^+\rightarrow\mathbb{N}$ be the function defined by: $w(\alpha_i+\cdots+\alpha_j)=(j-i+1)(n-j).$
Let $S=(\beta_1,\beta_2,\ldots,\beta_N)$ be a birational sequence where $\Phi^+=\{\beta_1,\ldots,\beta_N\}$. 
We consider the integral weight function 
$$
\Psi:\mathbb{Z}^N\rightarrow\mathbb{Z},\ \ \mathbf{m}\mapsto \sum_{i=1}^Nm_iw(\beta_i)
$$
and fix a lexicographic order ``$>$'' on $\mathbb{N}^N$. It is shown in \cite{FFR} that the global essential monoid 
$\Gamma(S,>)$ coincides with the monoid $\Gamma_\mathcal{P}$.
\par
In general, varying the function $w$ will change the monoid $\Gamma(S,>)$, yet in \cite{BFF} 
we give a family of such functions such that the global monoid stays constant, \emph{i.e.}, we have always $\Gamma(S,>)=\Gamma_\mathcal{P}$.
\item For various Weyl group elements $w\in W$ (\textit{triangular elements}), the Schubert variety $X_w$ degenerates to the irreducible toric variety defined by the $w$-face of $\mathcal{P}(\lambda)$ \cite{Fo2}.
This is similar to the case of $GT(\lambda)$, see Remark~\ref{rem-kempf}.
\item Let $G=SL_n$. We consider the following data:
\begin{itemize}
\item the birational sequence 
$$
S=(\alpha_1,\alpha_1+\alpha_2,\alpha_2,\alpha_1 + \alpha_2 +\alpha_3,\alpha_2 +\alpha_3,\alpha_3, \ldots,  \alpha_{n-1});
$$
\item an integral weight function $\Psi$ as in Example \ref{rootorder};
\item a $\Psi$-weighted lexicographic order on $\mathbb{N}^N$.
\end{itemize}
Then the corresponding global essential monoid $\Gamma(S,>)$ coincides with the monoid $\Gamma_\mathcal{P}$ (see \cite{FFL2}).
\end{enumerate}
 \end{remark}

\section{Gromov width of coadjoint orbits}\label{Gromov}
Let $\omega_{st}$ be the standard symplectic form on $\mathbb R^{2n}$. 
A famous theorem of Gromov (the {\it non-squeezing theorem}) affirms that a ball $B^{2n}(r)\subset (\mathbb R^{2n}, \omega_{st})$ 
cannot be symplectically embedded into $B^2(R)\times \mathbb R^{2n-2}\subset  (\mathbb R^{2n}, \omega_{st})$ unless $r\leq R$.
Given a symplectic manifold $M$, this result motivates the following questions:
\par\noindent
(1) what is the largest ball $B(r)$ that can be symplectically embedded into $M$? 
\par\noindent
(2) how many balls (of the same radius) can be symplectically packed into $M$? 

The first question is the source for the notation of the {\it Gromov width} of a $2n$-dimensional symplectic manifold $(M,\omega)$: it 
is the supremum of the set of $a$'s such that the ball of {\it capacity} $a$ (radius $\sqrt{\frac{a}{\pi}}$),
$$
B^{2n}_a = \big \{ (x_1,y_1,\ldots,x_n,y_n) \in  \mathbb R^{2n} \ \Big | \ \pi \sum_{i=1}^n (x_i^2+y_i^2) < a \big 
\} \subset  (\mathbb R^{2n}, \omega_{st}),
$$
 can be symplectically
embedded in $(M,\omega)$.  For the second question, a \emph{symplectic packing} of $(M,\omega)$ by $N$ equal balls of capacity $a$ is a symplectic embedding:
$$
\underbrace{B^{2n}_a\sqcup \cdots \sqcup B^{2n}_a}_{N\, {\rm times}}\rightarrow M
$$
of a disjoint union of $N$ balls of equal capacity $a$ into $M$. The Darboux theorem ensures that,
given $N$, such a packing always exists provided that the radius $r$ respectively the capacity $a$ is small enough. But when
one increases the capacity there will be obstructions for the existence of such a packing. A {\it full symplectic packing}
is a symplectic packing such that $N$ times the volume of $B^{2n}_a$ is equal to the volume of $(M,\omega)$.

An important class of symplectic manifolds is formed by the orbits of the coadjoint action of compact Lie groups.
Let $K$ be a compact Lie group, and  let $\mathfrak k^*$ be the dual of its Lie algebra $\mathfrak k$. 
Each orbit 
$\mathcal{O}\subset \mathfrak k^*$ of the coadjoint action of $K$ on $\mathfrak k^*$
is naturally equipped with the Kostant-Kirillov-Souriau symplectic form $\omega^{KKS}$, defined by:
$$\omega^{KKS}_{\xi}(X^\#,Y^\#)=\langle \xi, [X,Y]\rangle,\;\;\;\xi \in \mathcal{O} \subset \mathfrak k^*,\;X,Y \in \mathfrak k,$$
where $X^\#,Y^\#$ are the vector fields on $\mathfrak k^*$ induced by $X,Y \in \mathfrak k$ via the coadjoint action of $K$. Each coadjoint orbit passes through a point $\lambda$ in a positive Weyl chamber, we let $\mathcal{O}_\lambda$ denote it. Let $G=K_{\mathbb C}$ be the complexification of $K$. The coadjoint orbits of $K$ are diffeomorphic to $G/Q$
for some appropriate choices of a parabolic subgroup, and for $\lambda$ a dominant integral weight,
the pair $(\mathcal{O_\lambda}, \omega^{KKS})$ is symplectomorphic to $(G/Q,\omega^{FS})$, where the latter
denotes the Fubini-Study form induced by the embedding $G/Q\hookrightarrow \mathbb P(V(\lambda))$.

The momentum map for the standard $S=(S^1)^n$ action on $(\mathbb R^{2n}, \omega_{st})$ maps a ball of capacity $a$ into 
an $n$-dimensional simplex of size $a$, closed on $n$ sides: 
\begin{equation}\label{simplex}
\mathfrak S^n(a):=\{(x_1,\ldots,x_n) \in \mathbb R^n|\ 0\leq x_j< a,\  \sum_{j=1}^n x_j< a\}.
\end{equation}

Pabiniak has proved a kind of reverse implication of this fact, which provides a lower bound for the Gromov width
of a symplectic manifold: 
\begin{theorem}[\cite{Pa}]\label{embedding}
For any connected, proper (not necessarily compact) Hamiltonian $(S^1)^n$-space $M$ 
of dimension $2n$, with a momentum map $\Phi$, the Gromov width of $M$ is at least
$$\sup \{a>0\,|\, \exists \; \Psi \in GL(n,\mathbb Z), x \in
\mathbb R^n,\textrm{ such that }
\Psi (\textrm{int }\mathfrak S^n(a))+\,x \subset \Phi(M) \}.
$$
\end{theorem}

Of course, a coadjoint orbit (or $G/Q$ in the integral case) is not a toric variety. 
But combining this result with those of Harada and Kaveh \cite{HK} mentioned in the introduction,
Kaveh's result implies in the flag variety case (i.e. $\lambda$ is integral) the following result. To state it we need some notations: let $\Phi^+_Q=\{\beta_1,\ldots,\beta_N\}$ be the set of positive roots in the unipotent radical of $Q$, $S=(\beta_1,\beta_2,\ldots,\beta_N)$ be a good sequence (see section \ref{GoodOrdering}), $\Psi$ be the height function as in Example~\ref{rootorder} and take the right lexicographic ordering; the associated Newton-Okounkov body (see Definition \ref{Def:NOBody}) is denoted by $\Delta_\lambda (S)$.

\begin{theorem}[\cite{K2}]
The Gromov width of $(G/Q,\omega^{FS})$ is at least R, where R is the supremum of
the sizes of open simplices  that fit (up to $GL(N,\mathbb Z)$ transformation) in the interior of the Newton-Okounkov body $\Delta_\lambda (S)$, \emph{i.e.},
$$R= \sup \{r>0\,|\, \exists \; \Psi \in GL(N,\mathbb Z), x \in
\mathbb R^N,\textrm{ such that }\Psi (\textrm{int }\mathfrak S^N(r))+\,x \subset \Delta_\lambda (S) \}.$$
\end{theorem}

It is worthwhile to point out that only the Newton-Okounkov body is mentioned, there is no condition on the finite generation of the monoid associated to the valuation. We do not know in general whether the corresponding essential monoid $\Gamma(\lambda)$ is finitely generated, but one still has the equality with the valuation monoid
and one has a number of nice additional additive properties coming from the global essential monoid $\Gamma$, which provide the tools to prove:

\begin{theorem}[\cite{FLP}]\label{Prop:goodordering} Let $\lambda$ be an integral dominant weight.
For 
$$
k=\min\{\vert \langle\lambda,\alpha^\vee\rangle\vert \mid \alpha^{\vee} \textrm{ a  coroot and }\left\langle \lambda,\alpha^{\vee} \right\rangle \neq0\},
$$   
one has $\mathfrak S^N(k)\subset \Delta_\lambda$. In particular, the Gromov width of the coadjoint orbit through $\lambda$ is at least $k$.
\end{theorem}

With standard arguments one can extend this result to rational weights, the 
extension to irrational $\lambda$ is done by a ``Moser type" argument, described in detail in \cite{HP}. 
Combining this result with the upper bounds proved in \cite{CC}, one immediately obtains:

\begin{corollary}\label{corgw}
Let $K$ be a compact connected simple Lie group.
The Gromov width of a coadjoint orbit $\mathcal{O}_\lambda$ through $\lambda$, equipped with the Kostant-Kirillov-Souriau symplectic form, is given by 
\begin{equation}\label{gw formula}
\min\{\, \left|\left\langle \lambda,\alpha^{\vee} \right\rangle \right|\mid \  \alpha^{\vee} \textrm{ a  coroot and }\left\langle \lambda,\alpha^{\vee} \right\rangle \neq0\}.
 \end{equation}
\end{corollary}

Although this result had been proved in several cases (see \cite{FLP} for a detailed account), it was unsatisfactory that the proofs for the lower bounds were different for each group. The proof provided in \cite{FLP} is case-independent.
\par
In the same spirit, to construct symplectic packings can be reduced to study triangulations of the Newton-Okounkov body by simplices of a fixed size. In the following particular cases a unimodular triangulation of the body is known: (1). the Newton-Okounkov body is a stretched simplex (\cite{AKL, K2}); (2). the Newton-Okounkov body can be identified with some known polytopes such as the order polytope or chain polytope (\cite{Sta}). The first case can be applied to study the symplectic packing of the coadjoint orbit by balls of capacity no more than $1$ (\cite{K2}). We ask for a generalization of the second case to study the triangulations of marked order, marked chain or marked chain-order polytopes (\cite{FFMCO}), which can be then applied to the higher capacity cases.



\section{Small rank examples}\label{sec:smallrank}

\subsection{Symplectic case: \texorpdfstring{$\mathfrak{sp}_4$}{}}
We discuss briefly the case of toric degenerations of symplectic flag varieties in the rank $2$ case.
We consider three polytopes which
\begin{itemize}
\item are motivated from different setups;
\item give monomial bases of finite dimensional simple representations; \item lead to non-isomorphic toric degenerations of the symplectic flag variety associated to $\mathfrak{sp}_4$.
\end{itemize}
\par
Let $\lambda = m_1 \omega_1 + m_2 \omega_2$ be a dominant integral weight for $\mathfrak{sp}_4$.
\par
The first polytope is defined by:
\[
\mathbf{SP}_{4}(\lambda) = \left\{ (x_1,\ldots,x_4) \in \mathbb{R}_{\geq 0}^4 \mid 
{\small
\begin{array}{c} x_1\leq m_1, 2x_1+x_2+2x_3+2x_4\leq 2(m_1+m_2), \\  x_4\leq m_2, x_1+x_2+x_3+2x_4\leq m_1+2m_2\end{array}
}
\right\}.
\]
We denote the intersection $\mathbf{SP}_4(\lambda) \cap \mathbb{Z}^4$ by $\mathbf{SP}_4(\lambda) ^{\mathbb{Z}}$.
\begin{remark}\label{remark:sp4.1}
\begin{enumerate}
\item Up to permuting the second and the third coordinates, the polytope $\mathbf{SP}_4(\lambda)$ coincides with the one in Proposition 4.1 of \cite{K2}, 
which is unimodularly equivalent to the Newton-Okounkov body of a valuation arising from inclusions of (translated) Schubert varieties.
\item This polytope can be also constructed using any birational sequences $S=(\beta_1,\ldots,\beta_4)$ where $\Phi^+=\{\beta_1,\ldots,\beta_4\}$,
and the following $\Psi$-weighted function:
$$\Psi(\alpha_1) =1,\ \ \Psi(\alpha_2) = 1,\ \ \Psi(\alpha_1 + \alpha_2) = 2,\ \ \Psi(2\alpha_1 + \alpha_2) =1.$$
Let ``$>$'' be the lexicographic order. One obtains (\cite{BFF}) that the global essential monoid $\Gamma_1: =\Gamma(S,>)$ is finitely generated and saturated. 
\end{enumerate}
\end{remark}

The second polytope has been defined in \cite{FFL3} in the framework of PBW filtrations:
\[
\mathcal{Q}(\lambda) = \left\{ (x_1,\ldots,x_4) \in \mathbb{R}_{\geq 0}^4 \mid \begin{array}{c} x_1\leq m_1, x_1+x_2+x_3 \leq m_1 + m_2, 
\\  x_3\leq m_2, x_1+x_2+x_4\leq m_1+m_2.\end{array}\right\},
\]
denote again the intersection $\mathcal{Q}(\lambda) \cap \mathbb{Z}^4$ by $\mathcal{Q}_{\mathbb{Z}}(\lambda)$.
\begin{remark}\label{stringetal}
\begin{enumerate}
\item The very same polytope was suggested by \cite{ABS} in connection with the work on marked poset polytopes. 
It has been shown in \cite{Fo1} 
that this polytope is unimodularly equivalent to the generalized Gelfand-Tsetlin polytope, defined by Berenstein-Zelevinsky.
\item This polytope can be constructed using the birational sequence $S=(2\alpha_1+\alpha_2, \alpha_1+\alpha_2,\alpha_1,\alpha_2)$,
the integral weight function $\Psi$ as in Example \ref{homorder}, and, as total order, one takes the $\Psi$-weighted right opposite lexicographic order ``$>$''. Denote by $\Gamma_2:=\Gamma(S,>)$ the corresponding global essential monoid.
\item It has been further shown in \cite{L3} that the generalized Gelfand-Tsetlin polytope is unimodularly equivalent to the string polytope $Q_{s_2s_1s_2s_1}(\lambda)$.
\end{enumerate}
\end{remark}
The third polytope is the string polytope $Q_{s_1s_2s_1s_2}(\lambda)$, here the string cone is defined as
\[
C_{s_1s_2s_1s_2} = \left\{ (x_1,\ldots,x_4) \in \mathbb{R}^4_{\geq 0} \mid 2x_2 \geq x_3 \geq 2x_4 \right\}.
\]
Then the string polytope $Q_{s_1s_2s_1s_2}(\lambda)$ is defined to be the intersection of $C_{s_1s_2s_1s_2}$ with the half-spaces defined by the inequalities
\[
x_4 \leq m_2 \, , \, x_3 \leq m_1 + 2m_2 \, , \, x_2 \leq m_1  - m_2 \, , \, x_1 \leq 4m_1 - m_2.
\]

\begin{remark}\label{remark:sp4.3}
As shown in Section \ref{Sec:String}, these string polytopes can be realized using the birational sequence 
$S=(\alpha_1,\alpha_2,\alpha_1,\alpha_2)$, the height weighted function $\Psi$ as in Example \ref{rootorder}, and,
as total order, one takes the $\Psi$-weighted opposite lexicographic order $>$. 
Denote by $\Gamma_3:=\Gamma(S,>)$ the corresponding global essential monoid.

\end{remark}

It can be verified using POLYMAKE \cite{GJ}:
\begin{itemize}
\item The polytopes $\mathbf{SP}_4(\lambda)$, $\mathcal{Q}(\lambda)$, and $Q_{s_1s_2s_1s_2}(\lambda)$ are 
pairwise not unimodularly equivalent.
\item For $i=1,2,3$, $\Gamma_i$ (see Remark~\ref{remark:sp4.1}(2),\ref{stringetal}(2),\ref{remark:sp4.3}) is a finitely generated and saturated monoid. Let $X_i:=\text{Spec}(\mathbb{C}[\Gamma_i])$ denote the associated toric variety.
\item The toric varieties $X_1$, $X_2$ and $X_3$ are pairwise non-isomorphic.
\end{itemize}

\begin{remark}
By Remark~\ref{stringetal} we know that $\mathcal{Q}(\lambda)$, and $Q_{s_2s_1s_2s_1}(\lambda)$ are unimodularly equivalent to each other,
so the corollary states in particular that the two string polytopes $Q_{s_2s_1s_2s_1}(\lambda)$ and $Q_{s_1s_2s_1s_2}(\lambda)$ are not unimodularly equivalent.
\end{remark}

\subsection{The \texorpdfstring{$G_2$-case}{G2-case}}
Let $G$ be of type $G_2$. For a fixed dominant weight $\lambda=m_1\omega_1+m_2\omega_2\in\Lambda^+$, Gornitskii provided in \cite{G} a polytope $G(\lambda)$ defined by inequalities, whose lattice points parametrize monomial bases of the PBW-graded modules of simple modules. These polytopes $G(\lambda)$ can be realized as essential polytopes by considering
\begin{itemize}
\item the birational sequence $S=(3\alpha+2\beta, 3\alpha+\beta, 2\alpha+\beta, \alpha+\beta, \alpha, \beta)$ where $\alpha$ is the short simple root and $\beta$ is the long simple root;
\item the homogeneous integral weight function $\Psi$ as in Example \ref{homorder};
\item the $\Psi$-weighted lexicographic order ``$>$''.
\end{itemize}
Let $\Gamma:=\Gamma(S,>)$ be the associated global essential monoid. In this case, the monoid $\Gamma$ is finitely generated and saturated.
\par
We consider another setup consisting of
\begin{itemize}
\item the birational sequence $S'=(\alpha, 3\alpha+\beta, 2\alpha+\beta, 3\alpha+2\beta, \alpha+\beta, \beta)$;
\item the homogeneous integral weight function $\Psi':\mathbb{Z}^6\rightarrow\mathbb{Z}$ defined by:
$$
\Psi'(\mathbf m)=2m_1+m_2+3m_3+m_4+3m_5+m_6;
$$
\item the $\Psi'$-weighted lexicographic order ``$>'$''.
\end{itemize}
Let $\Gamma':=\Gamma(S',>')$ be the associated global essential monoid. The monoid $\Gamma'$ is not saturated \cite{BFF}.

\subsection{Toric degenerations of \texorpdfstring{$Gr_{2}(\mathbb C^4)$}{Gra\ss mannian of planes}}\label{gtwofour}

Let us examine how the results in previous sections can be applied to the special case of the Gra\ss mann variety $Gr_{2}(\mathbb C^4)$, 
giving rise to different toric degenerations. We keep the notations as in Section \ref{Sec:Grassmannian} for Gra{\ss}mann varieties.
\par
The homogeneous coordinate ring $\mathbb{C}[Gr_{2}(\mathbb C^4)]$ decomposes as $SL_4$-module:
$$\mathbb{C}[Gr_{2}(\mathbb C^4)]\cong \bigoplus_{k\geq 0}V(k\omega_2)^*.$$
As $SL_4$-modules, $V(\omega_2)\cong \Lambda^2 V(\omega_1)=\Lambda^2 \mathbb{C}^4$, and $e_1\wedge e_2$ is a highest weight vector.
\par
We consider the essential multi-exponents in $V(\omega_2)$ with respect to different ASM-sequences and show how they can be applied to the study of toric degenerations.

\subsubsection{The PBW polytope}
We fix the following data:
\begin{itemize}
\item a birational sequence ($\alpha_{1,3}$, $\alpha_{1,2}$, $\alpha_{2,3}$, $\alpha_1$, $\alpha_2$, $\alpha_3$) arising from a good ordering;
\item an integral weight function $\Psi$ as in Example \ref{homorder};
\item a $\Psi$-weighted right opposite lexicographic order on $\mathbb{N}^6$.
\end{itemize}

Let $es_1(\omega_2)$ denote the set of the corresponding essential multi-exponents in $V(\omega_2)$. By Example \ref{Ex:FFLSL4}, we have
$$es_1(\omega_2)=\{\bold{s}_1=(0,0,0,0,0,0),\ \bold{s}_2=(0,0,0,0,1,0),\ \bold{s}_3=(0,1,0,0,0,0),$$ $$\bold{s}_4=(0,0,1,0,0,0),\ \bold{s}_5=(1,0,0,0,0,0),\ \bold{s}_6=(1,0,0,0,1,0)\},$$
and the basis
$$\mathbb{B}_{\omega_2}=\{e_1\wedge e_2,\ e_1\wedge e_3=f_2\cdot e_1\wedge e_2,\ e_3\wedge e_2=f_{2,3}\cdot e_1\wedge e_2,$$
$$e_1\wedge e_4=f_{2,3}\cdot e_1\wedge e_2,\ e_4\wedge e_2=f_{1,3}\cdot e_1\wedge e_2,\ e_4\wedge e_3=f_{1,3}f_2\cdot e_1\wedge e_2\}.$$
Let $\Gamma_1(\omega_2)$ be the essential monoid associated to $\omega_2$. The dual essential basis elements are related to the Pl\"ucker coordinates by:
$$\bold{p_{[12]}}=\xi_{\omega_2,\bold{s}_1},\ \ \bold{p_{[13]}}=\xi_{\omega_2,\bold{s}_2},\ \ \bold{p_{[14]}}=\xi_{\omega_2,\bold{s}_4},$$
$$\bold{p_{[23]}}=-\xi_{\omega_2,\bold{s}_3},\ \ \bold{p_{[24]}}=-\xi_{\omega_2,\bold{s}_5},\ \ \bold{p_{[34]}}=-\xi_{\omega_2,\bold{s}_6}.$$
In the associated graded algebra, we have:
$$\xi_{\omega_2,\bold{s}_1}\cdot \xi_{\omega_2,\bold{s}_6}=\xi_{2\omega_2,(1,0,0,0,1,0)}=\xi_{\omega_2,\bold{s}_2}\cdot \xi_{\omega_2,\bold{s}_5},\ \  \xi_{\omega_2,\bold{s}_3}\cdot \xi_{\omega_2,\bold{s}_4}=\xi_{2\omega_2,(0,1,1,0,0,0)}.$$
This implies $\bold{p_{[12]}}\bold{p_{[34]}}-\bold{p_{[13]}}\bold{p_{[24]}}=0$, which defines the toric variety $\textrm{Spec}\,\mathbb{C}[\Gamma_1(\omega_2)]$. It is a toric degeneration of $Gr_{2}(\mathbb C^4)$.

\subsubsection{The string polytope}\label{Sec:Stringpol}
We fix the following data:
\begin{itemize}
\item a birational sequence ($\alpha_1$, $\alpha_2$, $\alpha_1$, $\alpha_3$, $\alpha_2$, $\alpha_1$);
\item an integral weight function $\Psi$ as in Example \ref{rootorder};
\item a $\Psi$-weighted opposite lexicographic order on $\mathbb{N}^6$.
\end{itemize}

Let $es_2(\omega_2)$ denote the set of the corresponding essential multi-exponents in $V(\omega_2)$:
$$es_2(\omega_2)=\{\bold{t}_1=(0,0,0,0,0,0),\ \bold{t}_2=(0,1,0,0,0,0),\ \bold{t}_3=(1,1,0,0,0,0),$$ $$\bold{t}_4=(0,0,0,1,1,0),\ \bold{t}_5=(1,0,0,1,1,0),\ \bold{t}_6=(0,1,1,1,1,0)\},$$
and the basis
$$\mathbb{B}_{\omega_2}=\{e_1\wedge e_2,\ e_1\wedge e_3=f_2\cdot e_1\wedge e_2,\ e_2\wedge e_3=f_1f_2\cdot e_1\wedge e_2,$$
$$e_1\wedge e_4=f_3f_2\cdot e_1\wedge e_2,\ e_2\wedge e_4=f_1f_3f_2\cdot e_1\wedge e_2,\ e_3\wedge e_4=f_2f_1f_3f_2\cdot e_1\wedge e_2\}.$$
Let $\Gamma_2(\omega_2)$ be the essential monoid associated to $\omega_2$. The dual essential basis elements are related to the Pl\"ucker coordinates by:
$$\bold{p_{[12]}}=\xi_{\omega_2,\bold{t}_1},\ \ \bold{p_{[13]}}=\xi_{\omega_2,\bold{t}_2},\ \ \bold{p_{[14]}}=\xi_{\omega_2,\bold{t}_4},$$
$$\bold{p_{[23]}}=\xi_{\omega_2,\bold{t}_3},\ \ \bold{p_{[24]}}=\xi_{\omega_2,\bold{t}_5},\ \ \bold{p_{[34]}}=\xi_{\omega_2,\bold{t}_6}.$$
In the associated graded algebra, we have:
$$\xi_{\omega_2,\bold{t}_2}\cdot \xi_{\omega_2,\bold{t}_5}=\xi_{2\omega_2,(1,1,0,1,1,0)}=\xi_{\omega_2,\bold{t}_3}\cdot \xi_{\omega_2,\bold{t}_4},\ \  \xi_{\omega_2,\bold{t}_1}\cdot \xi_{\omega_2,\bold{t}_6}=\xi_{2\omega_2,(0,1,1,1,1,0)}.$$
This implies $\bold{p_{[13]}}\bold{p_{[24]}}-\bold{p_{[14]}}\bold{p_{[23]}}=0$, which defines the toric variety $\textrm{Spec}\,\mathbb{C}[\Gamma_2(\omega_2)]$. It is a toric degeneration of $Gr_{2}(\mathbb C^4)$.

\subsubsection{The Lusztig polytope}
We fix the following data:
\begin{itemize}
\item a birational sequence ($\alpha_3$, $\alpha_{2,3}$, $\alpha_{1,3}$, $\alpha_2$, $\alpha_{1,2}$, $\alpha_1$) arising from the reduced decomposition $\underline{w}_0=s_3s_2s_1s_3s_2s_3$;
\item an integral weight function $\Psi$ as in Example \ref{rootorder};
\item a $\Psi$-weighted right opposite lexicographic order on $\mathbb{N}^6$.
\end{itemize}

Let $es_3(\omega_2)$ denote the set of the corresponding essential multi-exponents in $V(\omega_2)$:
$$es_3(\omega_2)=\{\bold{r}_1=(0,0,0,0,0,0),\ \bold{r}_2=(0,0,0,1,0,0),\ \bold{r}_3=(0,0,0,0,1,0),$$ $$\bold{r}_4=(1,0,0,1,0,0),\ \bold{r}_5=(1,0,0,0,1,0),\ \bold{r}_6=(0,1,0,0,1,0)\},$$
and the basis
$$\mathbb{B}_{\omega_2}=\{e_1\wedge e_2,\ e_1\wedge e_3=f_2\cdot e_1\wedge e_2,\ e_3\wedge e_2=f_{1,2}\cdot e_1\wedge e_2,$$
$$e_1\wedge e_4=f_3f_2\cdot e_1\wedge e_2,\ e_4\wedge e_2=f_3f_{1,2}\cdot e_1\wedge e_2,\ e_3\wedge e_4=f_{2,3}f_{1,2}\cdot e_1\wedge e_2\}.$$
Let $\Gamma_3(\omega_2)$ be the essential monoid associated to $\omega_2$. The dual essential basis elements are related to the Pl\"ucker coordinates by:
$$\bold{p_{[12]}}=\xi_{\omega_2,\bold{r}_1},\ \ \bold{p_{[13]}}=\xi_{\omega_2,\bold{r}_2},\ \ \bold{p_{[14]}}=\xi_{\omega_2,\bold{r}_4},$$
$$\bold{p_{[23]}}=-\xi_{\omega_2,\bold{r}_3},\ \ \bold{p_{[24]}}=-\xi_{\omega_2,\bold{r}_5},\ \ \bold{p_{[34]}}=\xi_{\omega_2,\bold{r}_6}.$$
In the associated graded algebra, we have:
$$\xi_{\omega_2,\bold{r}_2}\cdot \xi_{\omega_2,\bold{r}_5}=\xi_{2\omega_2,(1,0,0,1,1,0)}=\xi_{\omega_2,\bold{r}_3}\cdot \xi_{\omega_2,\bold{r}_4},\ \  \xi_{\omega_2,\bold{r}_1}\cdot \xi_{\omega_2,\bold{r}_6}=\xi_{2\omega_2,(0,1,0,0,1,0)}.$$
This implies $\bold{p_{[13]}}\bold{p_{[24]}}-\bold{p_{[14]}}\bold{p_{[23]}}=0$, which defines the toric variety $\textrm{Spec}\,\mathbb{C}[\Gamma_3(\omega_2)]$. It is a toric degeneration of $Gr_{2}(\mathbb C^4)$.

\subsubsection{The polytope arising from quantum PBW filtrations}
We fix the following data:
\begin{itemize}
\item a birational sequence ($\alpha_{1,3}$, $\alpha_2$, $\alpha_{2,3}$, $\alpha_1$, $\alpha_{1,2}$, $\alpha_3$) (in fact we can take an arbitrary birational sequence containing all positive roots);
\item an integral weight function $\Psi:\mathbb{Z}^6\rightarrow\mathbb{Z}$ determined by: for $\bold{m}=(m_1,\ldots,m_6)\in\mathbb{Z}^6$,
$$\Psi(\bold{m})=3m_1+2m_2+2m_3+3m_4+4m_5+m_6,$$
as in Remark \ref{remark-pbw}(5);
\item a $\Psi$-weighted lexicographic order on $\mathbb{N}^6$.
\end{itemize}

It is shown in \cite{FFR} that the essential basis associated to these data coincides with the one arising from the classical PBW filtration (Section \ref{GoodOrdering}).
\par
This integral weight function $\Psi$ endows degrees to the Pl\"ucker coordinates as follows:
$$\textrm{deg}(\bold{p_{[12]}})=\Psi((0,0,0,0,0,0))=0,\ \ \textrm{deg}(\bold{p_{[13]}})=\Psi((0,1,0,0,0,0))=2,$$
$$\textrm{deg}(\bold{p_{[14]}})=\Psi((0,0,1,0,0,0))=2,\ \ \textrm{deg}(\bold{p_{[23]}})=\Psi((0,0,0,0,1,0))=4,$$
$$\textrm{deg}(\bold{p_{[24]}})=\Psi((1,0,0,0,0,0))=3,\ \ \textrm{deg}(\bold{p_{[34]}})=\Psi((1,1,0,0,0,0))=5.$$
In the Pl\"ucker relation
$$
\mathbf p_{\mathbf{[14]}} \mathbf p_\mathbf{[23]} -
\mathbf p_{\mathbf{[13]}}\mathbf p_{\mathbf{[24]}}+
\mathbf p_{\mathbf{[12]}}\mathbf p_{\mathbf{[34]}}=0,
$$
$\textrm{deg}(\mathbf p_{\mathbf{[12]}} \mathbf p_\mathbf{[34]})=5=\textrm{deg}(\mathbf p_{\mathbf{[13]}} \mathbf p_\mathbf{[24]})$, but $\textrm{deg}(\mathbf p_{\mathbf{[14]}} \mathbf p_\mathbf{[23]})=6$. 

Hence in the associated graded algebra, $\mathbf p_{\mathbf{[12]}}\mathbf p_{\mathbf{[34]}}-\mathbf p_{\mathbf{[13]}}\mathbf p_{\mathbf{[24]}}=0$, 
defining the toric variety $\text{Spec\,}\mathbb{C}[\Gamma_1(\omega_2)]$. It is a toric degeneration of $Gr_{2}(\mathbb C^4)$.

\subsubsection{Remark}
We have seen four different approaches to degenerate $Gr_{2}(\mathbb C^4)$ into toric varieties. 
Due to the fact that the dimension is rather small, it turns out that all these approaches are giving isomorphic toric degenerations. 
This can be easily verified as the defined polytopes are unimodularly equivalent, which does not hold in general.
\par
For example, for the Gra\ss mann variety $Gr_{3}(\mathbb C^6)$, the toric degenerations arising from the PBW polytope described in Section \ref{GoodOrdering} and from the string polytope in Section \ref{Sec:String} associated to the reduced decomposition
$$\underline{w}_0=s_1s_2s_1s_3s_2s_1s_4s_3s_2s_1s_5s_4s_3s_2s_1$$
are not isomorphic since the corresponding polytopes are not unimodularly equivalent \cite{Fo1}. More on the Gra\ss mann variety $Gr_{3}(\mathbb C^6)$ can be found in Example~\ref{GR36}.

\section{Toric degenerations and cluster varieties}\label{cluster}
 A new and  different approach towards toric degenerations has been suggested recently by 
Gross, Hacking, Keel and Kontsevich \cite{GHKK} in the case of cluster algebras which satisfy certain additional properties. 
The authors believe that all the elementary constructions of toric geometry extend to log Calabi-Yau varieties (with
maximal boundary), and they prove many such results in the case of cluster varieties. 
As a guiding example for the connection with flag varieties and their likes, one should have in mind $\mathbb C[SL_n/U]$.
By Gei\ss, Leclerc and Schr\"oer \cite{GLS},  it admits a cluster algebra structure, and by a recent 
preprint of Magee \cite{M}, has the desired additional properties. 
To go into all the technical details would blow up the framework of this overview, so we will only give simplified versions 
of the statements and restrict ourselves to some special cases. Especially we will only consider cluster algebras of geometric type, \emph{i.e.}, with the exchange matrix being skew-symmetric. For more details on cluster algebras see also \cite{FZ,GLS,Z}.

\subsection{Cluster duality and mirror symmetry} \label{RietschWilliams}  
To see the connection between the ways polytopes show up in Section~\ref{Section:birational} and how polytopes
show up in the cluster context, let us start with a recent article by Rietsch and Williams \cite{RW}.
The authors do not discuss degenerations, but their approach is inspiring to see how cluster algebras and toric charts 
can help to understand the connection between two very different ways of looking at polytopes.

Rietsch and Williams investigate the case of $X=Gr_{n-k}(\mathbb C^n)$, the Gra\ss mann variety of $(n-k)$-dimensional 
subspaces of $\mathbb C^n$. They consider coordinate charts on $X$ and the ``mirror dual'' Gra\ss mann variety 
$\check{X}:=Gr_{k}\left(\left(\mathbb{C}^n\right)^*\right)$, respectively. The charts on both sides are 
obtained from a choice of combinatorial objects, a reduced plabic graph $\mathcal{G}$\index{$\mathcal{G}$, plabic graph} with trip permutation $\pi_{k,n}$. 

The coordinate system on $X$ is given by an injective map
\begin{equation}\label{eq:ytransformation}
\Phi_{\mathcal G}\index{$\Phi_{\mathcal G}$, the positive chart associated to a plabic graph $\mathcal{G}$}: \left(\mathbb{C}^*\right)^{\mathcal{P}_{\mathcal{G}}}\rightarrow Gr_{n-k}(\mathbb C^n)
\end{equation} 
constructed in \cite{T}. Here $\mathcal{P}_{\mathcal{G}}$\index{$\mathcal{P}_{\mathcal{G}}$, Pl\"ucker coordinates associated to plabic graph $\mathcal{G}$} is an index set for a certain set of Pl\"ucker coordinates read off 
from the graph $\mathcal{G}$ by a combinatorial rule. The image of the restriction of $\Phi_{\mathcal G}$ to 
$\left(\mathbb{R}_{>0}\right)^{\mathcal{P}_{\mathcal{G}}}$ is the totally positive Gra\ss mann variety in its Pl\"ucker 
embedding and thus Rietsch and Williams refer to it as a {\it positive chart}.

To each positive chart $\Phi_{\mathcal G}$ and $r>0$, the authors associate
a \textit{Newton-Okounkov-type polytope} $NO^r(\mathcal G)$. They define (similar as in Section~\ref{NObody}) $NO^1(\mathcal G)$
as the convex hull of certain lattice points which correspond to the vanishing behavior of sections of the line bundle
$\mathcal L$, the ample generator of ${\rm Pic\,} X$. For any $r>1$, $NO^r(\mathcal{G})$ is defined to be the $r$-fold Minkowski sum of $NO^1(\mathcal{G})$.\footnote{The polytope $NO^r(\mathcal{G})$ is contained in the Newton-Okounkov body associated to the line bundle $\mathcal{L}^r$, but may be different from it.}
 
It is known \cite{Sc} that the homogeneous coordinate ring of (the affine cone over) $\check{X}$ admits a cluster algebra structure. 
Each cluster ${\bf x}=(x_1,x_2,\ldots,x_m)$\index{$\bold{x}$, a cluster} of this cluster algebra gives rise to a toric chart on $\check{X}$
\begin{equation}\label{eq:toricchart}
\check{X}_{\bf x}\index{$\check{X}_{\bf x}$, toric chart associated to a cluster $\bold{x}$}:=\{y\in \check{X} \mid x_i(y)\ne 0, \ 1 \le i \le m\}.
\end{equation}
The index set $\mathcal{P}_{\mathcal{G}}$ labels a collection of Pl\"ucker coordinates that form an entire cluster (recall that in general there will be many clusters 
where the variables do not consist entirely of Pl\"ucker coordinates). Therefore, from (\ref{eq:toricchart}), we get a map
\begin{equation}\label{eq:xtransformation}
{\Phi_{\mathcal G}^{\vee}}: \left(\mathbb{C}^*\right)^{\mathcal{P}_{\mathcal{G}}} \rightarrow \check{X}
\end{equation}
called {\it cluster chart} which satisfies $p_{\nu}(\Phi_{\mathcal G}^{\vee}\left((t_{\mu})_{\mu})\right)=t_{\nu}$ for 
$\nu\in \mathcal{P}_{\mathcal{G}}$ and $p_{\nu}$ the associated Pl\"ucker coordinate.

In mirror symmetry the {\it mirror} of the Gra\ss mann variety $X$ is a Landau-Ginzburg 
model, which can be described as the pair $(\check{X}^o,W_q)$, where $\check X^o$ 
is the complement of a particular anticanonical divisor in the Langlands dual Gra\ss mann variety $\check{X}$, 
and $W_q$ is a regular map on $\check{X}^o$, called superpotential. The image of $\Phi_{\mathcal G}^{\vee}$ lies in $\check{X}^o$. The condition of the tropicalized version of the superpotential $W_{t^r}\circ \Phi_{\mathcal G}^{\vee}$ (\emph{i.e.}, the superpotential written in a cluster expansion in terms of the cluster consisting of Pl\"ucker coordinates on $\check{X}$ labeled by $\mathcal{P}_{\mathcal{G}}$ and replacing the $q$-variable by $t^r$) to have non-negative value gives rise to a set of linear inequalities defining a polytope $Q_{\mathcal G}^r$.

This is the important difference between the two constructions of the polytopes.
In the second procedure the 
description of the polytope is given as the intersection of half-spaces while in the first procedure the description of the polytope is given by taking a convex hull of a set of lattice points.
The main result in \cite{RW} is:
\begin{theorem}
The two polytopes $NO^r(\mathcal G)$ and $Q_{\mathcal G}^r$ coincide  
for all reduced plabic graphs $\mathcal G$ with trip permutation $\pi_{k,n}$ and all $r>0$. 
\end{theorem}
More precisely, for the respective charts (\ref{eq:ytransformation}) and (\ref{eq:xtransformation}) associated to a particular reduced plabic graph $\mathcal G_0$ Rietsch and Williams prove that $NO^r(\mathcal G_0)$ and $Q_{\mathcal G_0}^r$
coincide. Moreover, they claim that one can show that the polytopes are unimodularly equivalent to 
Gelfand-Tsetlin polytopes for $r\omega_{n-k}$,  where $\omega_{n-k}$ is the $SL_n$-fundamental weight corresponding to the Gra\ss mann variety $X$. 
To switch from one chart to another, one has a combinatorial procedure which describes how to transform one reduced plabic graph with trip permutation $\pi_{k,n}$ to another. This induces a piecewise linear map describing how to switch from one associated polytope to the 
other, called a {\it tropicalized cluster mutation}. They show that if $\mathcal G'$ and $\mathcal G''$ are plabic graphs
which are related by a single move, then the piecewise-linear transformation relating 
$NO^r(\mathcal G')$ to $NO^r(\mathcal G'')$ is the
same as the piecewise-linear transformation relating $Q_{\mathcal G'}^r$ to $Q_{\mathcal G''}^r$. 
It follows that for any reduced plabic graph one has $NO^r(\mathcal G')=Q_{\mathcal G'}^r$.
\begin{example}\rm\label{GR36}
We consider the Gra\ss mann variety $Gr_{3}(\mathbb C^6)$. In \cite{RW} it is shown
that the Newton-Okounkov-type polytope associated to the particular reduced plabic graph $\mathcal G_0$
is unimodularly equivalent to the Gelfand-Tsetlin polytope $GT(\omega_3)$. By performing a 
tropicalized cluster mutation on the unique square in this particular plabic graph 
(see Section 6 and 7 in \cite{RW}), one obtains a Newton-Okounkov-type polytope 
which is unimodularly equivalent to the polytope given by the following 15 constraints:
\[
{\tiny
T(\omega_3) := \left\{ (x_i) \in \mathbb{R}^9 \mid \begin{array}{l}  -x_9 + x_5 \leq 1 \; ; \; x_2 - x_5 - x_6 + x_9 \geq 0\; ; \; x_4 - x_5 - x_8 + x_9 \geq 0 \\ 
-x_4 - x_7 + x_8 \geq 0\; ; \; -x_2 - x_3 + x_6 \geq 0    \; ; \;  x_1 - x_4 \geq 0    \; ; \; x_1 - x_2 \geq 0\\  
-x_2 + x_3 \geq 0\; ; \; -x_4 + x_7 \geq 0   \; ; \; -x_1 + x_5 \geq 0 \\ 
-x_1 + x_2 + x_4 \geq 0   \; ; \; -x_1 + x_2 - x_5 + x_8 \geq 0    -x_1 + x_4 - x_5 + x_6 \geq 0  \\
 -x_1 - x_2 + x_6 \geq 0     \; ; \; -x_1 - x_4 + x_8 \geq 0   
 \end{array}
\right\}.
}
\]
This polytope is normal. It is neither unimodularly equivalent to $GT(\omega_3)$ nor to $\mathcal{P}(\omega_3)$ (Section \ref{GoodOrdering}). The polytopes $T(\omega_3)$ and $\mathcal{P}(\omega_3)$ even share the same $f$-vector. Up to a unimodularly equivalence, this polytope $T(\omega_3)$ can be realized as an essential polytope associated to some birational sequence and appropriate total order.
\par
Applying mutations on the reduced plabic graphs (\cite{RW}) gives rise to 34 different reduced plabic graphs, and the associated Newton-Okounkov-type polytopes have 6 isomorphism classes. Five of these classes appear in the framework of tropical Gra\ss mannians \cite{SS} and can be also realized using the construction of birational sequences.
\end{example}

\begin{example}\label{FFLVpolytope}\rm
We consider the Gra\ss mann variety $Gr_{k}(\mathbb C^n)$ and the particular reduced plabic graph $\mathcal G_0$. The Newton-Okounkov-type polytope associated to this plabic graph is claimed to be unimodularly equivalent to the Gelfand-Tsetlin polytope $GT(\omega_k)$ for $SL_n$ (\cite{RW}). This plabic graph corresponds to a seed $\bold{s}$ with cluster $\{\bold{p}_{J_1},\ldots,\bold{p}_{J_s}\}$ in the cluster algebra structure on $Gr_k(\mathbb{C}^n)$ (\cite{Sc}), where $\bold{p}_{J_r}$ is the Pl\"ucker coordinate associated to $J_r\in I_{k,n}$. Let $J_r^\vee\in I_{n-k,n}$ be the complement of $J_r$ in $\{1,2,\ldots,n\}$. The set of Pl\"ucker coordinates $\{\bold{p}_{J_1^\vee},\ldots,\bold{p}_{J_s^\vee}\}$ satisfies the maximally weakly separated condition \cite{OPS}, therefore it is the cluster in a seed $\bold{s}^\vee$ in the cluster algebra structure on $Gr_{n-k}(\mathbb{C}^n)$, hence corresponds to a plabic graph $\mathcal G_0^m$.
\par
It is proved in \cite{FF} (see Conjecture 1 in \cite{Fang}) that the Newton-Okounkov-type polytope obtained from the plabic graph $\mathcal{G}_0^m$ is unimodularly equivalent to the polytope $\mathcal{P}(\omega_{n-k})$ (see Section \ref{GoodOrdering}) for $SL_n$. This provides a duality between the Gelfand-Tsetlin polytopes and polytopes arising from PBW filtration in the case of Gra\ss mannian. 
\end{example}

\subsection{Cluster varieties, toric charts and Langlands dual}\label{LLdual}
The basic objects in \cite{GHKK} are cluster ensembles introduced in \cite{FG} and related spaces. To define these varieties, the definition of a cluster algebra is reformulated, not including the cluster variables in the seed datum but rather defining it coordinate free. 
Recall that we assume for simplicity our cluster algebra to be of geometric type.  A {\it seed datum} consists of
\begin{itemize}
\item a lattice $N$ with a skew-symmetric bilinear form $\{\cdot, \cdot\}: N \times N \rightarrow \mathbb{Q}$;
\item a saturated {\it unfrozen sublattice} $N_{\rm uf}\subset N$;
\item an index set $I$ with $|I|=\text{rank}N$, $I_{\rm uf}\subset I$ with $|I_{\rm uf}|=\text{rank}N_{\rm uf}$;
\item the dual lattice $M=\Hom(N,\mathbb{Z})$ to $N$;
\item a basis ${\bf s}:=(e_i \mid i \in I)$ of $N$ such that $\{e_i \mid i \in I_{\rm uf}\}$ is a basis for $N_{\rm uf}$.
\end{itemize}
For a choice of seed ${\bf s}$ we further obtain a dual basis $\{{e_i^*}\}$ for $M$.
To a seed datum ${\bf s}$ we associate two tori in the spirit of (\ref{eq:toricchart}):
$$\begin{matrix}
{\mathcal X}_{\bf s}=T_M=\text{Spec\,}\mathbb{C}[N] &  \text{ and } & {\mathcal A}_{\bf s}=T_N=\text{Spec\,}\mathbb{C}[M].
\end{matrix}$$
For a seed data we can define a mutation of the seed data in direction $i$ for all $i \in I_{\rm uf}$. The $(e_i \mid i \in I)$ change according to the usual exchange relation given by the matrix of skew-symmetric form with respect to that choice of basis $\{\cdot, \cdot\}$ (see \cite{FZ2}, Equation 2.4) while the ${e_i^{*}}$ change according to the $Y$-seed mutation (see \cite{FZ2}, Equation 2.3). This induces birational maps on the associated tori. We define the $\mathcal X$\index{$\mathcal X$, $\mathcal X$-cluster variety} and $\mathcal A$\index{$\mathcal A$, $\mathcal A$-cluster variety}-cluster variety by gluing the tori obtained by all choices of seeds using these birational maps. The pair ($\mathcal X$, $\mathcal A$) is called a {\it cluster ensemble} in \cite{FG}. We define
\begin{equation}\label{dualFGvariety}
{\mathcal A}_\Gamma^\vee :={\mathcal X}_{\Gamma}  \quad 
{\mathcal X}_\Gamma^\vee :={\mathcal A}_{\Gamma}.
\end{equation}

\begin{example}\rm In the case of the Gra\ss mann variety the positive charts (\ref{eq:ytransformation}) in Section \ref{RietschWilliams} are instances of charts for 
an $\mathcal X$-cluster variety and the cluster charts (\ref{eq:xtransformation}) are instances of charts for an $\mathcal A$-cluster variety.
\end{example}

The {\it $\mathcal{X}$-cluster algebra} ({\it $\mathcal A$-cluster algebra}) associated with a seed ${\bf s}$ is defined to be $\Gamma(\mathcal X, \mathcal{O}_{\mathcal X})$ ($\Gamma(\mathcal A, \mathcal{O}_{\mathcal A})$).  The $\mathcal A$-cluster algebra is usually called the upper cluster algebra (see \cite{BFZ}).

Often the coordinate ring of
a variety $X$ has a {\it cluster algebra} structure, but $X$ is not a cluster variety in the sense above.
In these cases one tries to establish $X$ as a {\it partial compactification} $\overline{\mathcal A}$ of a cluster variety.
An example of such an instance is $SL_n/U$, which can be thought of as a partial compactification of the double Bruhat cell
$G^{e,w_0}$ (see Example~\ref{doublebruhat}).

\begin{example}\label{doublebruhat}\rm 
Let $G=SL_n(\mathbb C)$ with Borel subgroup $B^+$ the subgroup of upper triangular matrices,
and let $B^-$ be the opposite Borel subgroup of lower triangular matrices. Let $U$ be the unipotent radical
of $B^+$. For $u, v$ in the Weyl group $W$ of $G$, the double Bruhat cell
$G^{u,v}$ is defined to be $G^{u,v} = (B^+uB^+) \cap (B^-vB^-)$. It has been proved by
Berenstein, Fomin and Zelevinsky \cite{BFZ} that the coordinate ring of any double Bruhat cell in a
semisimple complex Lie group is naturally isomorphic to an upper cluster algebra. Further the double Bruhat cell $G^{e,w_0}$ is a cluster variety.
There is an open embedding of $G^{e,w_0}$ into $G/U$ given by $g\mapsto g^t U$. 
One can view $G/U$ (up to codimension 2 differences, but this does not affect the corresponding coordinate rings) 
as a partial compactification of $G^{e,w_0}$, afforded by taking frozen variables to $0$. So the coordinate ring of $G/U$
inherits in a natural way a cluster algebra structure from $G^{e,w_0}$.
\end{example}
Let $\mathbb Z^T$\index{$\mathbb Z^T$, tropical semifield $\mathbb Z$} denote the tropical semifield $\mathbb Z$ with  operations ``$\max$'' and ``$+$'',
and for a cluster variety $\mathcal A$ let $\mathcal A(\mathbb Z^T)$ be the $\mathbb Z^T$-valued points 
(also called the tropical points).
Whenever one fixes a seed $\mathbf s$, then one gets an 
identification of $\mathcal A(\mathbb Z^T)$ and $\mathcal X(\mathbb Z^T)$ with the corresponding
underlying lattice $N$ respectively $M$ (see \cite{GHKK}).
Nevertheless the tropicalized versions of the birational mutation maps induce only piecewise linear maps between 
the various identifications, not linear maps. The idea behind the approach of \cite{FG}, \cite{GHKK} is
that despite this many of the elementary constructions in toric geometry still work. For example,
the coordinate ring of $T_{N}$ has a canonical basis, the characters of the torus, which are canonically
indexed by the character lattice $M$. Correspondingly, the algebra of regular functions on $\mathcal A$
should come with a canonical basis indexed by the tropical points in $\mathcal X(\mathbb Z^T)$. 
This is the Fock-Goncharov-conjecture in a weak form.
Similarly, for every fixed seed one has an identification of $\mathcal X(\mathbb Z^T)$ with the underlying 
lattice $M$.                     
Again, the identification of $\mathcal X (\mathbb Z^T)$ with the lattices associated to the various charts 
induce only piecewise linear (and not linear) identifications.

\subsection{A canonical basis and polytopes}\label{canbasis}
Although the Fock-Goncharov-con\-jecture is false in general (see \cite{GHK}), in certain particularly nice cases, the cluster algebra comes equipped with such a canonical basis as mentioned above. 
Conjecturally the construction of a canonical basis in \cite{GHKK} holds for a large class of cluster algebras of representation theoretical interest, for example for 
all $G/U$, $G$ a semisimple algebraic group. We would like to emphasize that this construction is about cluster algebras, 
the construction of the canonical basis in \cite{GHKK} works without any representation theory.
The question whether the construction above provides a basis or not is related to the question whether 
the full Fock-Goncharov-conjecture holds for the cluster variety respectively its partial compactification.
 
At the moment it is known to hold for $SL_n/U$ (see  \cite{M}). 
In this case the cluster variety is the double Bruhat cell $\mathcal A = G^{e,w_0}$
and we consider $G/U$ as a partial compactification as in the example above. 
In this case the coordinate ring of the cluster variety $\mathcal A$ comes equipped with a canonical basis $\mathbb B$ constructed in \cite{GHKK}.
The elements of the basis are naturally indexed by the $\mathbb Z^T$-valued points (the tropical points) of the mirror cluster variety,
which is in this case the cluster variety $\mathcal X$. Now viewing $SL_n/U$ as a partial  compactification
of the cluster variety $\mathcal A$, this endows the mirror cluster variety $\mathcal X$ with a canonical potential $W$. The tropicalization $W^T:\mathcal X(\mathbb R^T) \rightarrow \mathbb R$ is a piecewise linear map, and the  condition to be non-negative: $\{x\in \mathcal X(\mathbb R^T) \mid W^T\ge 0\}$
cuts out a tropical cone. This means that for any choice of a seed and the associated identification of $\mathcal X(\mathbb R^T)$ with a real vector space 
and  $\mathcal X(\mathbb Z^T)$ with a lattice, $\{x\in \mathcal X(\mathbb R^T) \mid W^T\ge 0\}$ is a (usual) cone in a real vector space. A beautiful
consequence of the theory is:
\begin{theorem}\cite{GHKK,M}
The integer points in the tropical cone defined by the condition $ W^T\ge 0$ parametrize a basis $\Theta$ of the coordinate 
ring of the partial compactification $\overline{\mathcal A}=SL_n/U$. 
\end{theorem}
This description resembles that of the string cones (see section~\ref{Sec:String}), whose lattice points parametrize
the canonical basis constructed by Kashiwara and Lusztig. 
Indeed, Magee has shown in \cite{M} that for a particular choice of a seed (and the corresponding identification of 
$\mathcal X(\mathbb R^T)$ with a real vector space) one obtains a cone 
which is unimodularly equivalent to the Gelfand-Tsetlin cone.
Gross, Hacking, Keel and Kontsevich conjecture that all string cones can be obtained in this way.

We expect that among the cones which one gets via the construction of Gross, Hacking, Keel and Kontsevich
one will find many cones which are not unimodularly equivalent to string cones and which are interesting for 
representation theoretic considerations.

Let $D$ be the fixed maximal torus in $SL_n$.  The canonical basis $\Theta$ consists of eigenvectors with respect to the right 
action of $D$ on $SL_n/U$. Given a dominant integral weight $\lambda$, the irreducible representation $V(\lambda)\subset \mathbb C[SL_n/U]$
is a $D$-eigenspace for this right action, and the intersection $\Theta_\lambda=\Theta\cap V(\lambda)$ is a canonical basis
of this representation space.
On the combinatorial side, one can associate to $\lambda$ an affine tropical subspace of $\mathcal X(\mathbb Z^T)$,
such that the intersection with the tropical cone defines a (tropical) polytope $P(\lambda)$. 
\begin{theorem}\cite{GHKK,M}
The lattice points of the polytope $P(\lambda)$
parametrize the elements of the canonical basis $\Theta_\lambda\subset V(\lambda)$. 
Again, for a particular seed and the identification
$\mathcal X(\mathbb Z^T)=M$ above, the polytope $P(\lambda)\subset M_{\mathbb R}$ is 
unimodularly equivalent to the corresponding Gelfand-Tsetlin pattern
associated to the dominant weight $\lambda$.
\end{theorem}
\subsection{From polytopes and cluster varieties to degenerations}\label{polyclusterdegen}
Let $\mathcal A$ be a cluster variety. An important tool used by  Gross, Hacking, Keel and Kontsevich
is the principle cluster algebra $\mathcal A_{prin}$. In terms of the notation introduced in section~\ref{LLdual},
this procedure glues together the $\mathcal{A}$- and the $\mathcal{X}$-cluster varieties considering the new variables coming from the lattice $N$ as frozen variables (\emph{i.e.}, treating them as principal coefficients). For details and the notation see \cite{GHK}. The dual $\mathcal X_{prin}$ is obtained the same way now taking the new variables coming from the lattice $M$ as frozen variables.

A choice of an initial seed $\mathbf s$ provides a partial compactification $\mathcal A_{prin,\mathbf s}$
of $\mathcal A_{prin}$ by allowing the principal coefficients to be zero.
One gets an induced flat map $\pi:\mathcal A_{prin,\mathbf s}\rightarrow \mathbb A^n$
with $\mathcal A$ being the fiber over $(1,\ldots,1)$, and the central fiber $\pi^{-1}(0) \subset \mathcal A_{prin,\mathbf s}$ 
is the algebraic torus $T_{N}$ given by the seed $\mathbf s$.

%

Each regular function $f:{\mathcal X}_{prin}\rightarrow \mathbb A^1$ has a canonical piecewise linear
tropicalisation $f^T:{\mathcal X}_{prin}(\mathbb R^T)\rightarrow \mathbb R$, which is conjectured to be convex.
Roughly speaking, this means that for any seed $\mathbf s$ and the corresponding identification of
${\mathcal X}_{prin}(\mathbb R^T)$ with a real vector space, $f^T$ can be described as the minimum function
for a finite number of linear functions, $f^T$ satisfies a certain convexity condition, and a condition of the differential
of the piecewise linear function is satisfied. These important properties of convex piecewise linear functions imply:
\begin{equation}\label{clusterpolytope}
\Xi_f:=\{x\in \mathcal A^\vee_{prin}(\mathbb R^T)\mid f^T(x)\ge -1\}
\end{equation}
is a convex polytope.

Assume now $\Xi_f$ is bounded, rational, and satisfies an additional positivity condition. One has a natural projection map 
$\rho : \mathcal X_{prin}\rightarrow \mathcal X$  with tropicalization
$\rho^T:\mathcal X_{prin}(\mathbb Z^T ) \rightarrow \mathcal X(\mathbb Z^T )$. Let
$\overline{\Xi}_f\subset \mathcal X(\mathbb R^T)$ be the image of $\Xi_f$ with respect to $\rho^T$, then $\overline{\Xi}_f$
is a bounded polytope with $0$ in its interior.

With the help of the polytope $\Xi_f$ and using a Rees-type construction for graded rings, Gross, Hacking, Keel and Kontsevich construct a graded ring $\tilde R$ 
together with a flat morphism 
\begin{equation}\label{flatfamily}
{\rm Proj\,}(\tilde R)\rightarrow \mathbb A^n,
\end{equation}
such that:
\begin{theorem}\cite{GHKK}\label{GHKKfibration}
The central fiber of this map is the polarized toric variety for the torus $T_{N}$ given by the polyhedron $\overline{\Xi}_f$
and the generic fiber is a compactification of the cluster variety $\mathcal A$.
\end{theorem}
The polytope can be chosen so that the boundary of the compactification is very
simple, a union of toric varieties \cite{GHKK}. As an example, the authors mention the
open subset $Gr^o_k(\mathbb C^n)\subset Gr_k(\mathbb C^n)$
which is a cluster variety. Then generic compactifications given by bounded
polytopes give an alternative compactification of this open subset in which one replaces certain Schubert
cells (which are highly non-toric) by toric varieties.

Gross, Hacking, Keel and Kontsevich \cite{GHKK} conjecture that all the constructions of toric degenerations of flag varieties 
(and there likes) in \cite{AB} (see also section~\ref{ALBriondeg}) are instances of the theorem above. 
This would imply that for a given string polytope one can find an $f$ such that the generic
fiber above is a generalized flag variety $G/P_\lambda$ and $\overline \Xi_f$ is unimodularly equivalent to the given string polytope.

\section{An outlook: conjectures and questions}\label{openquestion}
We have presented several different methods to construct flat degenerations of a flag variety into a toric variety.
It is natural to ask for connections between them.
\subsection{Cluster varieties and birational sequences}
Having a birational sequence \eqref{birational1} $\pi:Z_S\rightarrow U^-$ for $U^-$, the natural embedding
of $U^-$ into $G/B$ gives hence a map $\pi:Z_S\rightarrow G/B$. As a variety, $Z_S$ is just an affine space $\mathbb A^N$
containing a torus $\mathbb T^N$. In the case where it is known that $G/B$ admits an open cluster subvariety, it is natural
to ask: under which conditions is the image $\pi(\mathbb T^N)$ a toric chart
for the cluster subvariety of $G/B$? If this holds, how does the choice of a weighted lexicographic order  (Definition~\ref{weightorder})
and hence the construction of the global essential monoid and the associated Newton-Okounkov body
correspond to the choice of a regular function on the cluster subvariety of $G/B$ (as in section~\ref{polyclusterdegen})
and its associated polytope \eqref{clusterpolytope}? 
A first result in this direction has been obtained by Magee \cite{M} in the case of $SL_n(\mathbb C)$
for a particular seed. Roughly speaking, this problem is also related to the 
question pursued by Rietsch and Williams in the case of the Gra\ss mann variety: what is the connection between the polytopes 
obtained by the Newton-Okounkov approach and the polytopes obtained using the superpotential $W_q$ (\cite{RW}, see 
also section~\ref{RietschWilliams}). 

\subsection{Properties of birational sequences}
Given a connected complex reductive group $G$, it is natural to ask
for a parametrization of all possible birational sequences. 
The next obvious question is to ask for a characterization of those
weighted lexicographic orderings associated to such a sequence
such that the essential monoid is finitely generated respectively finitely
generated and saturated. We conjecture that the essential monoid is
always finitely generated, independently of the choice of the birational sequence
and the choice of the ordering.

For applications it might be useful to have available all possible choices
of birational sequences and weighted orders. Is it possible to characterize 
a special class of birational sequences and orderings (like, for example, the PBW-case
and the right opposite lexicographic order?) which gives all $T$-equivariant toric degenerations?
Is it possible to characterize the polytopes arising from toric degenerations of flag varieties, or, 
more generally, arising from toric degenerations of spherical varieties?

Schubert varieties and Richardson varieties are important geometric objects in the representation theory closely related to flag varieties,
it would be interesting to obtain a characterization of those birational sequences 
which are compatible with such varieties. 

\subsection{Generalizations}
Is it possible to realize Chriv\`i's degeneration (Section \ref{Section:Chiviri}) in a framework similar to that of 
birational sequences (section \ref{Section:birational})? 
Is it possible to extend the methods using birational sequences to  affine Kac-Moody groups \cite{FeFiR}?
Is it possible to generalize the filtration aspect to the setting of quantum groups \cite{FFR}? 
Some recent work in this direction can be found in \cite{BFF}. As both have a Gr\"obner
background, it might be interesting to study the relation between the quantum degree cone \cite{BFF} and
the tropical Gra\ss mannian \cite{SS}, or, in general, tropical flag varieties. 

Instead of studying toric degenerations, one may also look at intermediate steps,
as, for example, the PBW-degenerate flag variety in \cite{CL, FeFiL,F1}.  
A related class of degenerations are the linear degenerations of flag varieties studied in \cite{CFFFR},
which are degenerations of the flag variety as a quiver Gra\ss mannian. In connection with the cluster variety
approach it is natural to ask whether it is possible to identify these degenerations as a fiber of the flat morphism 
constructed by Gross, Hacking, Keel and Kontsevich (see \eqref{flatfamily} and Theorem~\ref{GHKKfibration}).

\printindex

\end{document}